%

%
\ifx\begin\undefined\else\global\advance\srcdepth by
1\expandafter\endinput\fi

\def\begin{}
\newcount\srcdepth
\srcdepth=1

\outer\def\bye{\global\advance\srcdepth by -1
  \ifnum\srcdepth=0
    \def\endcmd{\vfill\supereject\nopagenumbers\par\vfill\supereject\end}
  \else\def\endcmd{}\fi
  \endcmd
}



\def\initialize#1#2#3#4#5#6{
  \ifnum\srcdepth=1
  \magnification=#1
  \hsize = #2
  \vsize = #3
  \hoffset=#4
  \advance\hoffset by -\hsize
  \divide\hoffset by 2
  \advance\hoffset by -1truein
  \voffset=#5
  \advance\voffset by -\vsize
  \divide\voffset by 2
  \advance\voffset by -1truein
  \advance\voffset by #6
  \baselineskip=13pt
  \emergencystretch = 0.05\hsize
  \fi
}

\def\print{\initialize{1095}
  {5.5truein}{8.5truein}{8.5truein}{11truein}{-.0625truein}}

\newif\ifblackboardbold

\blackboardboldtrue


\font\titlefont=cmbx12 scaled\magstephalf
\font\sectionfont=cmbx12

\font\scriptit=cmti10 at 7pt
\font\scriptsl=cmsl10 at 7pt
\scriptfont\itfam=\scriptit
\scriptfont\slfam=\scriptsl


\newfam\bboldfam
\ifblackboardbold
\font\tenbbold=msbm10
\font\sevenbbold=msbm7
\font\fivebbold=msbm5
\textfont\bboldfam=\tenbbold
\scriptfont\bboldfam=\sevenbbold
\scriptscriptfont\bboldfam=\fivebbold
\def\bbold{\fam\bboldfam\tenbbold}
\else
\def\bbold{\bf}
\fi


\newfam\msamfam
\font\tenmsam=msam10
\font\sevenmsam=msam7
\font\fivemsam=msam5
\textfont\msamfam=\tenmsam
\scriptfont\msamfam=\sevenmsam
\scriptscriptfont\msamfam=\fivemsam

\newfam\msbmfam
\font\tenmsbm=msam10
\font\sevenmsbm=msam7
\font\fivemsbm=msam5
\textfont\msbmfam=\tenmsbm
\scriptfont\msbmfam=\sevenmsbm
\scriptscriptfont\msbmfam=\fivemsbm

\newcount\amsfamcount 
\newcount\classcount   
\newcount\positioncount
\newcount\codecount
\newcount\n             
\def\newsymbol#1#2#3#4#5{               
\n="#2                                  
\ifnum\n=1 \amsfamcount=\msamfam\else   
\ifnum\n=2 \amsfamcount=\msbmfam\else   
\ifnum\n=3 \amsfamcount=\eufmfam
\fi\fi\fi
\multiply\amsfamcount by "100           
\classcount="#3                 
\multiply\classcount by "1000           
\positioncount="#4#5            
\codecount=\classcount                  
\advance\codecount by \amsfamcount      
\advance\codecount by \positioncount
\mathchardef#1=\codecount}              


\font\Arm=cmr9
\font\Ai=cmmi9
\font\Asy=cmsy9
\font\Abf=cmbx9
\font\Brm=cmr7
\font\Bi=cmmi7
\font\Bsy=cmsy7
\font\Bbf=cmbx7
\font\Crm=cmr6
\font\Ci=cmmi6
\font\Csy=cmsy6
\font\Cbf=cmbx6

\ifblackboardbold
\font\Abbold=msbm10 at 9pt
\font\Bbbold=msbm7
\font\Cbbold=msbm5 at 6pt
\fi

\def\small{%
\textfont0=\Arm \scriptfont0=\Brm \scriptscriptfont0=\Crm
\textfont1=\Ai \scriptfont1=\Bi \scriptscriptfont1=\Ci
\textfont2=\Asy \scriptfont2=\Bsy \scriptscriptfont2=\Csy
\textfont\bffam=\Abf \scriptfont\bffam=\Bbf \scriptscriptfont\bffam=\Cbf
\def\rm{\fam0\Arm}\def\mit{\fam1}\def\oldstyle{\fam1\Ai}%
\def\bf{\fam\bffam\Abf}%
\ifblackboardbold
\textfont\bboldfam=\Abbold
\scriptfont\bboldfam=\Bbbold
\scriptscriptfont\bboldfam=\Cbbold
\def\bbold{\fam\bboldfam\Abbold}%
\fi
\rm
}








\newlinechar=`@
\def\forwardmsg#1#2#3{\immediate\write16{@*!*!*!* forward reference should
be: @\noexpand\forward{#1}{#2}{#3}@}}
\def\nodefmsg#1{\immediate\write16{@*!*!*!* #1 is an undefined reference@}}

\def\forwardsub#1#2{\def\newref{{#2}{#1}}}

\def\forward#1#2#3{%
\expandafter\expandafter\expandafter\forwardsub\expandafter{#3}{#2}
\expandafter\ifx\csname#1\endcsname\relax\else%
\expandafter\ifx\csname#1\endcsname\newref\else%
\forwardmsg{#1}{#2}{#3}\fi\fi%
\expandafter\let\csname#1\endcsname\newref}

\def\firstarg#1{\expandafter\argone #1}\def\argone#1#2{#1}
\def\secondarg#1{\expandafter\argtwo #1}\def\argtwo#1#2{#2}

\def\ref#1{\expandafter\ifx\csname#1\endcsname\relax
  {\nodefmsg{#1}\bf`#1'}\else
  \expandafter\firstarg\csname#1\endcsname
  ~\expandafter\secondarg\csname#1\endcsname\fi}

\def\refs#1{\expandafter\ifx\csname#1\endcsname\relax
  {\nodefmsg{#1}\bf`#1'}\else
  \expandafter\firstarg\csname #1\endcsname
  s~\expandafter\secondarg\csname#1\endcsname\fi}

\def\refn#1{\expandafter\ifx\csname#1\endcsname\relax
  {\nodefmsg{#1}\bf`#1'}\else
  \expandafter\secondarg\csname #1\endcsname\fi}



\def\widow#1{\vskip 0pt plus#1\vsize\goodbreak\vskip 0pt plus-#1\vsize}



\def\marginlabel#1{}

\def\showlabelsabove{
\font\labelfont=cmss10 at 6pt
\def\marginlabel##1{\rlap{\smash{\raise 10pt\hbox{\labelfont##1}}}}
}

\newcount\seccount
\newcount\proccount
\seccount=0
\proccount=0

\def\stdskip{\vskip 9pt plus3pt minus 3pt}
\def\stdbreak{\par\removelastskip\penalty-100\stdskip}

\def\proof{\stdbreak\noindent{\bf Proof. }}

\def\qed{\vrule height 1.2ex width .9ex depth .1ex}

\def\Box{
  \ifmmode\eqno\qed
  \else\ifvmode\removelastskip\line{\hfil\qed}
  \else\unskip\quad\hskip-\hsize
    \hbox{}\hskip\hsize minus 1em\qed\par
  \fi\stdbreak\fi}

\def\references{
  \removelastskip
  \widow{.05}
  \vskip 24pt plus 6pt minus 6 pt
  \leftline{\sectionfont References}
  \nobreak\stdskip\noindent}

\def\ifempty#1#2\endB{\ifx#1\endA}
\def\makeref#1#2#3{\ifempty#1\endA\endB\else\forward{#1}{#2}{#3}\fi}

\outer\def\section#1 #2\par{
  \removelastskip
  \global\advance\seccount by 1
  \global\proccount=0\relax
                \edef\numtoks{\number\seccount}
  \makeref{#1}{Section}{\numtoks}
  \widow{.05}
  \vskip 24pt plus 6pt minus 6 pt
  \message{#2}
  \leftline{\marginlabel{#1}\sectionfont\numtoks\quad #2}
  \nobreak\stdskip}

\def\proclamation#1#2{
  \outer\expandafter\def\csname#1\endcsname##1 ##2\par{
  \stdbreak
  \global\advance\proccount by 1
  \edef\numtoks{\number\seccount.\number\proccount}
  \makeref{##1}{#2}{\numtoks}
  \noindent{\marginlabel{##1}\bf #2 \numtoks\enspace}
  {\sl##2\par}
  \stdbreak}}

\def\othernumbered#1#2{
  \outer\expandafter\def\csname#1\endcsname##1{
  \stdbreak
  \global\advance\proccount by 1
  \edef\numtoks{\number\seccount.\number\proccount}
  \makeref{##1}{#2}{\numtoks}
  \noindent{\marginlabel{##1}\bf #2 \numtoks\enspace}}}

\proclamation{definition}{Definition}
\proclamation{lemma}{Lemma}
\proclamation{proposition}{Proposition}
\proclamation{theorem}{Theorem}
\proclamation{corollary}{Corollary}
\proclamation{conjecture}{Conjecture}

\othernumbered{example}{Example}
\othernumbered{remark}{Remark}
\othernumbered{construction}{Construction}


\def\figure#1{
 \global\advance\figcount by 1
 \goodbreak
 \midinsert#1\smallskip
 \centerline{Figure~\number\figcount}
 \endinsert}

\def\capfigure#1#2{
 \global\advance\figcount by 1
 \goodbreak
 \midinsert#1\smallskip
 \vbox{\small\noindent {\bf Figure~\number\figcount:} #2}
 \endinsert}

\def\capfigurepair#1#2#3#4{
 \goodbreak
 \midinsert
 #1\smallskip
 \global\advance\figcount by 1
 \vbox{\small\noindent {\bf Figure~\number\figcount:} #2}
 \vskip 12pt
 #3\smallskip
 \global\advance\figcount by 1
 \vbox{\small\noindent {\bf Figure~\number\figcount:} #4}
 \endinsert}


\def\baretable#1#2{
\vbox{\offinterlineskip\halign{
 \strut\kern #1\hfil##\kern #1
 &&\kern #1\hfil##\kern #1\cr
 #2
}}}

\def\gridtablesub#1#2#3{
\vbox{\offinterlineskip\halign{
 \strut\vrule\kern #1\hfil##\hfil\kern #2\vrule
 &&\kern #1\hfil##\kern #2\vrule\cr
 \noalign{\hrule}
 #3
 \noalign{\hrule}
}}}






\newif\iftextures
\input epsf

\newcount\figcount
\figcount=0
\newcount\figxscale
\newcount\figyscale
\newcount\figxoffset
\newcount\figyoffset
\newbox\drawing
\newcount\drawbp
\newdimen\drawx
\newdimen\drawy
\newdimen\ngap
\newdimen\sgap
\newdimen\wgap
\newdimen\egap

\def\drawbox#1#2#3{\vbox{
  \epsfgetbb{#2.eps} 
  \drawbp=\epsfurx
  \advance\drawbp by-\epsfllx\relax
  \multiply\drawbp by #1
  \divide\drawbp by 100
  \drawx=\drawbp bp
  \drawbp=\epsfury
  \advance\drawbp by-\epsflly\relax
  \multiply\drawbp by #1
  \divide\drawbp by 100
  \drawy=\drawbp bp
  \iftextures
  		\figxscale=#1
    \multiply\figxscale by 10
    \setbox\drawing=\vbox to \drawy{\vfil
      \hbox to \drawx{\special{illustration #2.eps scaled
\number\figxscale}\hfil}}
  \else 
    \figxoffset=-\epsfllx
    \multiply\figxoffset by#1
    \divide\figxoffset by100
    \figyoffset=-\epsflly
    \multiply\figyoffset by#1
    \divide\figyoffset by100
    \setbox\drawing=\vbox to \drawy{\vfil
      \hbox to \drawx{\includegraphics{#2.eps}\hfil}}
  \fi
  \setbox\drawing=\vbox{\offinterlineskip\box\drawing\kern 0pt}
   \ngap=0pt \sgap=0pt \wgap=0pt \egap=0pt
  \setbox0=\vbox{\offinterlineskip
    \box\drawing \ifgridlines\drawgrid\drawx\drawy\fi #3}
  \kern\ngap\hbox{\kern\wgap\box0\kern\egap}\kern\sgap}}

\def\draw#1#2#3{
  \setbox\drawing=\drawbox{#1}{#2}{#3}
  \global\advance\figcount by 1
  \edef\numtoks{\number\figcount}
  \makeref{fig:#2}{Figure}{\numtoks}
  \goodbreak
  \midinsert
  \centerline{\ifgridlines\boxgrid\drawing\fi\box\drawing}
  \smallskip
  \vbox{\offinterlineskip
    \centerline{Figure~\number\figcount}
    \smash{\marginlabel{#2}}}
  \endinsert}

\def\capdraw#1#2#3#4{
  \setbox\drawing=\drawbox{#1}{#2}{#3}
  \global\advance\figcount by 1
  \edef\numtoks{\number\figcount}
  \makeref{fig:#2}{Figure}{\numtoks}
  \goodbreak
  \midinsert
  \centerline{\ifgridlines\boxgrid\drawing\fi\box\drawing}
  \smallskip
  \vbox{\offinterlineskip
    \vskip 4pt
    \vbox{\centerline{\lineskip=3pt\small\noindent
                       {\bf Figure~\number\figcount:} #4}}
    \smash{\marginlabel{fig:#2}}}
  \endinsert}

\def\capdrawpair#1#2#3#4#5#6#7#8{
  \goodbreak
  \midinsert
  \setbox\drawing=\drawbox{#1}{#2}{#3}
  \global\advance\figcount by 1
  \edef\numtoks{\number\figcount}
  \makeref{fig:#2}{Figure}{\numtoks}
  \centerline{\ifgridlines\boxgrid\drawing\fi\box\drawing}
  \smallskip
  \vbox{\offinterlineskip
    \vskip 4pt
    \vbox{\lineskip=3pt\small\noindent {\bf Figure~\number\figcount:} #4}
    \smash{\marginlabel{fig:#2}}}
  \vskip 12pt
  \setbox\drawing=\drawbox{#5}{#6}{#7}
  \global\advance\figcount by 1
  \edef\numtoks{\number\figcount}
  \makeref{fig:#6}{Figure}{\numtoks}
  \centerline{\ifgridlines\boxgrid\drawing\fi\box\drawing}
  \smallskip
  \vbox{\offinterlineskip
    \vskip 4pt
    \vbox{\lineskip=3pt\small\noindent {\bf Figure~\number\figcount:} #8}
    \smash{\marginlabel{fig:#6}}}
  \endinsert}

\def\nextfigtoks{%
  \advance\figcount by 1%
  \edef\numtoks{\number\figcount}%
  \advance\figcount by -1}

\def\nextfig{\nextfigtoks Figure~\numtoks}

\newif\ifgridlines
\newbox\figtbox
\newbox\figgbox
\newdimen\figtx
\newdimen\figty

\newdimen\bwd
\bwd=2sp 

\def\hline#1{\vbox{\smash{\hbox to #1{\leaders\hrule height \bwd\hfil}}}}

\def\vline#1{\hbox to 0pt{%
  \hss\vbox to #1{\leaders\vrule width \bwd\vfil}\hss}}

\def\clap#1{\hbox to 0pt{\hss#1\hss}}
\def\vclap#1{\vbox to 0pt{\offinterlineskip\vss#1\vss}}

\def\hstutter#1#2{\hbox{%
  \setbox0=\hbox{#1}%
  \hbox to #2\wd0{\leaders\box0\hfil}}}

\def\vstutter#1#2{\vbox{
  \setbox0=\vbox{\offinterlineskip #1}
  \dp0=0pt
  \vbox to #2\ht0{\leaders\box0\vfil}}}

\def\crosshairs#1#2{
  \dimen1=.002\drawx
  \dimen2=.002\drawy
  \ifdim\dimen1<\dimen2\dimen3\dimen1\else\dimen3\dimen2\fi
  \setbox1=\vclap{\vline{2\dimen3}}
  \setbox2=\clap{\hline{2\dimen3}}
  \setbox3=\hstutter{\kern\dimen1\box1}{4}
  \setbox4=\vstutter{\kern\dimen2\box2}{4}
  \setbox1=\vclap{\vline{4\dimen3}}
  \setbox2=\clap{\hline{4\dimen3}}
  \setbox5=\clap{\copy1\hstutter{\box3\kern\dimen1\box1}{6}}
  \setbox6=\vclap{\copy2\vstutter{\box4\kern\dimen2\box2}{6}}
  \setbox1=\vbox{\offinterlineskip\box5\box6}
  \smash{\vbox to #2{\hbox to #1{\hss\box1}\vss}}}

\def\boxgrid#1{\rlap{\vbox{\offinterlineskip
  \setbox0=\hline{\wd#1}
  \setbox1=\vline{\ht#1}
  \smash{\vbox to \ht#1{\offinterlineskip\copy0\vfil\box0}}
  \smash{\vbox{\hbox to \wd#1{\copy1\hfil\box1}}}}}}

\def\drawgrid#1#2{\vbox{\offinterlineskip
  \dimen0=\drawx
  \dimen1=\drawy
  \divide\dimen0 by 10
  \divide\dimen1 by 10
  \setbox0=\hline\drawx
  \setbox1=\vline\drawy
  \smash{\vbox{\offinterlineskip
    \copy0\vstutter{\kern\dimen1\box0}{10}}}
  \smash{\hbox{\copy1\hstutter{\kern\dimen0\box1}{10}}}}}

\def\figtext#1#2#3#4#5{
  \setbox\figtbox=\vbox{\hbox{#5}\kern 0pt}
  \figtx=-#3\wd\figtbox \figty=-#4\ht\figtbox
  \advance\figtx by #1\drawx \advance\figty by #2\drawy
  \dimen0=\figtx \advance\dimen0 by\wd\figtbox \advance\dimen0 by-\drawx
  \ifdim\dimen0>\egap\global\egap=\dimen0\fi
  \dimen0=\figty \advance\dimen0 by\ht\figtbox \advance\dimen0 by-\drawy
  \ifdim\dimen0>\ngap\global\ngap=\dimen0\fi
  \dimen0=-\figtx
  \ifdim\dimen0>\wgap\global\wgap=\dimen0\fi
  \dimen0=-\figty
  \ifdim\dimen0>\sgap\global\sgap=\dimen0\fi
  \smash{\rlap{\vbox{\offinterlineskip
    \hbox{\hbox to \figtx{}\ifgridlines\boxgrid\figtbox\fi\box\figtbox}
    \vbox to \figty{}
    \ifgridlines\crosshairs{#1\drawx}{#2\drawy}\fi
    \kern 0pt}}}}

\def\swtext#1#2#3{\figtext{#1}{#2}00{#3}}
\def\setext#1#2#3{\figtext{#1}{#2}10{#3}}

\def\wtext#1#2#3{\figtext{#1}{#2}0{.5}{#3}}
\def\etext#1#2#3{\figtext{#1}{#2}1{.5}{#3}}
\def\ntext#1#2#3{\figtext{#1}{#2}{.5}1{#3}}
\def\stext#1#2#3{\figtext{#1}{#2}{.5}0{#3}}


\def\hpad#1#2#3{\hbox{\kern #1\hbox{#3}\kern #2}}
\def\vpad#1#2#3{\setbox0=\hbox{#3}\vbox{\kern #1\box0\kern #2}}

\def\wpad#1#2{\hpad{#1}{0pt}{#2}}
\def\epad#1#2{\hpad{0pt}{#1}{#2}}
\def\npad#1#2{\vpad{#1}{0pt}{#2}}
\def\spad#1#2{\vpad{0pt}{#1}{#2}}

\def\swpad#1#2#3{\spad{#1}{\wpad{#2}{#3}}}
\def\sepad#1#2#3{\spad{#1}{\epad{#2}{#3}}}




\def\stack#1#2#3{\vbox{\offinterlineskip
  \setbox2=\hbox{#2}
  \setbox3=\hbox{#3}
  \dimen0=\ifdim\wd2>\wd3\wd2\else\wd3\fi
  \hbox to \dimen0{\hss\box2\hss}
  \kern #1
  \hbox to \dimen0{\hss\box3\hss}}}


\def\hexp#1{%
  \setbox0=\hbox{${}^{#1}$}%
  \hbox to .5\wd0{\box0\hss}}

\def\hsub#1{%
  \setbox0=\hbox{${}_{#1}$}%
  \hbox to .5\wd0{\box0\hss}}



\def\bmatrix#1#2{{\left(\vcenter{\halign
  {&\kern#1\hfil$##\mathstrut$\kern#1\cr#2}}\right)}}

\def\rightarrowmat#1#2#3{
  \setbox1=\hbox{\small\kern#2$\bmatrix{#1}{#3}$\kern#2}
  \,\vbox{\offinterlineskip\hbox to\wd1{\hfil\copy1\hfil}
    \kern 3pt\hbox to\wd1{\rightarrowfill}}\,}

\def\leftarrowmat#1#2#3{
  \setbox1=\hbox{\small\kern#2$\bmatrix{#1}{#3}$\kern#2}
  \,\vbox{\offinterlineskip\hbox to\wd1{\hfil\copy1\hfil}
    \kern 3pt\hbox to\wd1{\leftarrowfill}}\,}

\def\rightarrowbox#1#2{
  \setbox1=\hbox{\kern#1\hbox{\small #2}\kern#1}
  \,\vbox{\offinterlineskip\hbox to\wd1{\hfil\copy1\hfil}
    \kern 3pt\hbox to\wd1{\rightarrowfill}}\,}

\def\leftarrowbox#1#2{
  \setbox1=\hbox{\kern#1\hbox{\small #2}\kern#1}
  \,\vbox{\offinterlineskip\hbox to\wd1{\hfil\copy1\hfil}
    \kern 3pt\hbox to\wd1{\leftarrowfill}}\,}








\def\quiremacro#1#2#3#4#5#6#7#8#9{
  \expandafter\def\csname#1\endcsname##1{
  \ifnum\srcdepth=1
  \magnification=#2
  \input quire
  \hsize=#3
  \vsize=#4
  \htotal=#5
  \vtotal=#6
  \shstaplewidth=#7
  \shstaplelength=#8
  \hoffset=\htotal
  \advance\hoffset by -\hsize
  \divide\hoffset by 2
  \ifnum\vsize<\vtotal
    \voffset=\vtotal
    \advance\voffset by -\vsize
    \divide\voffset by 2
  \fi
  \advance\voffset by #9
  \shhtotal=2\htotal
  \baselineskip=13pt
  \emergencystretch = 0.05\hsize
  \horigin=0.0truein
  \vorigin=0.0truein
  \shthickness=0pt
  \shoutline=0pt
  \shcrop=0pt
  \shvoffset=-1.0truein
  \ifnum##1>0\quire{#1}\else\qtwopages\fi
  \fi
}}



\quiremacro{letterbooklet} 
{1000}{4.79452truein}{7truein}{5.5truein}{8.5truein}{0.01pt}{0.66truein}
{-.0625truein}

\quiremacro{Afourbooklet}
{1095}{5.25truein}{7truein}{421truept}{595truept}{0.01pt}{0.66truein}
{-.0625truein}

\quiremacro{legalbooklet}
{1095}{5.25truein}{7truein}{7.0truein}{8.5truein}{0.01pt}{0.66truein}
{-.0625truein}

\quiremacro{twoupsub} 
{895}{4.5truein}{7truein}{5.5truein}{8.5truein}{0pt}{0pt}{.0625truein}


\quiremacro{Afourviewsub} 
{1000}{5.0228311in}{7.7625571in}{421truept}{595truept}{0.1pt}{0.5\vtotal}
{-.0625truein}


\quiremacro{viewsub}
{1095}{5.5truein}{8.5truein}{461truept}{666truept}{0.1pt}{0.5\vtotal}
{-.125truein}


\newcount\countA
\newcount\countB
\newcount\countC

\def\monthname{\begingroup
  \ifcase\number\month
    \or January\or February\or March\or April\or May\or June\or
    July\or August\or September\or October\or November\or December\fi
\endgroup}

\def\dayname{\begingroup
  \countA=\number\day
  \countB=\number\year
  \advance\countA by 0 
  \advance\countA by \ifcase\month\or
    0\or 31\or 59\or 90\or 120\or 151\or
    181\or 212\or 243\or 273\or 304\or 334\fi
  \advance\countB by -1995
  \multiply\countB by 365
  \advance\countA by \countB
  \countB=\countA
  \divide\countB by 7
  \multiply\countB by 7
  \advance\countA by -\countB
  \advance\countA by 1
  \ifcase\countA\or Sunday\or Monday\or Tuesday\or Wednesday\or
    Thursday\or Friday\or Saturday\fi
\endgroup}

\def\timename{\begingroup
   \countA = \time
   \divide\countA by 60
   \countB = \countA
   \countC = \time
   \multiply\countA by 60
   \advance\countC by -\countA
   \ifnum\countC<10\toks1={0}\else\toks1={}\fi
   \ifnum\countB<12 \toks0={\sevenrm AM}
     \else\toks0={\sevenrm PM}\advance\countB by -12\fi
   \relax\ifnum\countB=0\countB=12\fi
   \hbox{\the\countB:\the\toks1 \the\countC \thinspace \the\toks0}
\endgroup}

\def\timestamp{\dayname, \the\day\ \monthname\ \the\year, \timename}


\print



\def\COMMENT#1\par{\bigskip\hrule\smallskip#1\smallskip\hrule\bigskip}

\def\enma#1{{\ifmmode#1\else$#1$\fi}}

\def\mathbb#1{{\bbold #1}}
\def\mathbf#1{{\bf #1}}


\def\NN{\enma{\mathbb{N}}}

\def\PP{\enma{\mathbb{P}}}
\def\RR{\enma{\mathbb{R}}}
\def\ZZ{\enma{\mathbb{Z}}}


\def\cEE{\enma{\cal E}}
\def\cHH{\enma{\cal H}}
\def\cFF{\enma{\cal F}}
\def\cII{\enma{\cal I}}
\def\cOO{\enma{\cal O}}

\def\aa{\enma{\mathbf{a}}}
\def\bb{\enma{\mathbf{b}}}
\def\cc{\enma{\mathbf{c}}}
\def\dd{\enma{\mathbf{d}}}
\def\ee{\enma{\mathbf{e}}}

\def\pp{\enma{\mathbf{p}}}
\def\qq{\enma{\mathbf{q}}}

\def\uu{\enma{\mathbf{u}}}
\def\vv{\enma{\mathbf{v}}}

\def\xx{\enma{\mathbf{x}}}
\def\yy{\enma{\mathbf{y}}}

\def\DD{\enma{\mathbf{D}}}
\def\FF{\enma{\mathbf{F}}}

\def\zero{\enma{\mathbf{0}}}


\def\Cl{\enma{\rm{Cl}}}

\def\set#1{\enma{\{#1\}}}
\def\setdef#1#2{\enma{\{\;#1\;\,|\allowbreak
  \;\,#2\;\}}}

\def\idealdef#1#2{\enma{\langle\;#1\;\,
  |\;\,#2\;\rangle}}

\def\mtext#1{\;\,\allowbreak\hbox{#1}\allowbreak\;\,}
\def\abs#1{\enma{\left| #1 \right|}}
\def\floor#1{\enma{\lfloor  #1 \rfloor}}
\def\ceil#1{\enma{\lceil  #1 \rceil}}


\def\codim{\mathop{\rm cod}\nolimits}
\def\hull{\mathop{\rm hull}\nolimits}

\def\im{\mathop{\rm im}\nolimits}

\def\rank{\mathop{\rm rank}\nolimits}

\def\Ext{\mathop{\rm Ext}\nolimits}
\def\Coh{\mathop{\rm Coh}\nolimits}
\def\ini{\mathop{\rm in}\nolimits}
\def\fib{\mathop{\rm fib}\nolimits}

\newsymbol\boxtimes1202


\forward{introd}{Section}{1}
\forward{unimodularity}{Section}{2}
\forward{main}{Section}{3}
\forward{mono}{Section}{4}
\forward{examples}{Section}{5}
\forward{diagonal}{Section}{6}
\forward{beilinson}{Section}{7}
\forward{fig:initial}{Figure}{4}

\overfullrule=0pt

\hbox{}
\smallskip
\centerline{\titlefont Syzygies of Unimodular Lawrence Ideals} \bigskip
\centerline{Dave Bayer \quad Sorin Popescu \quad Bernd Sturmfels}

\bigskip
\bigskip

\noindent {Abstract:}
{\sl Infinite hyperplane arrangements whose vertices form a lattice
are studied from the point of view of commutative algebra. The quotient
of such an arrangement modulo the lattice action represents the minimal
free resolution of the associated binomial ideal, which defines
a toric subvariety in a product of projective lines.
Connections to graphic arrangements and
to Beilinson's spectral sequence are explored.}

\bigskip

\section{introd} {Introduction}

\noindent  We are interested in the defining ideals of
toric subvarieties in a product of projective lines $\, \PP^1 \times
\PP^1 \times \cdots \times \PP^1$.  Writing $(x_i:y_i)$ for the
homogeneous coordinates of the $i$-th factor $\PP^1$, these are the
binomial ideals in $2n$ variables of the following form:
$$ J_L \quad = \quad \idealdef{\xx^\aa \yy^\bb - \xx^\bb \yy^\aa}{\aa-\bb
\in L} \quad \subset \quad
S \, = \,  k[x_1,\ldots,x_n, \, y_1,\ldots,y_n], $$
where $L$ is a sublattice of $\ZZ^n$ and $k$ is a field.  Here
$\,\xx^\aa = \,x_1^{a_1} x_2^{a_2} \cdots x_n^{a_n} \,$
for $\aa = (a_1,\ldots,a_n) \in \NN^n$. Binomial ideals
of the form $J_L$ are called {\it Lawrence ideals}.
They provide the algebraic analogue to the
Lawrence construction for convex polytopes [Zi, \S 6.6].
Lawrence polytopes enjoy remarkable
rigidity properties, such as [Zi, Theorem 6.27]. On the algebraic
side, rigidity of Lawrence ideals manifests itself in the
following result, which appears in [Stu, Theorem 7.1].
Recall for part (d) that the {\it Graver basis} consists of all 
binomials $\xx^\aa \yy^\bb - \xx^\bb \yy^\aa$ in $J_L$
such that the only vector $\aa' - \bb' \in L \backslash \set{\zero}$ 
with $\zero \leq \aa'\leq \aa$ and $\zero \leq \bb' \leq \bb$ is 
$\aa - \bb$ itself.

\proposition{lawrencegens} The following sets of binomials in a
Lawrence ideal $J_L$ coincide:
\item{\rm (a)} Any minimal set of binomial generators of $J_L$.
\item{\rm (b)} Any reduced Gr\"obner basis for $J_L$.
\item{\rm (c)} The universal Gr\"obner basis for $J_L$ (the union of all reduced
Gr\"obner bases).
\item{\rm (d)} The Graver basis for $J_L$.

Cellular resolutions, as defined in [BS], provide a natural geometric
framework for studying homological, algorithmic and combinatorial
properties of monomial and binomial ideals. One instance is the 
{\it Scarf complex}, which gives the minimal resolution for generic
monomial ideals [BPS] and generic lattice ideals [PS]. 
The {\it hull resolution} of [BS] generalizes the Scarf complex
and provides a cellular resolution for 
arbitrary co-Artinian monomial modules; however, it need not be
minimal. For lattice modules, the
 hull resolution is compatible with the lattice action and determines a 
cellular resolution of the corresponding lattice ideal [BS, Theorem 3.9].

In this paper we present a minimal cellular resolution,
which happens to also coincide with the hull resolution, for
a remarkable class of nongeneric lattice ideals.
These are the unimodular Lawrence ideals $J_L$, which are characterized as follows:

\theorem{unimod} For a sublattice $L$ of $\ZZ^n$ the following
conditions are equivalent:
\item{\rm (a)} The Lawrence ideal $J_L$ possesses an initial monomial
  ideal which is radical.
\item{\rm (b)} Every initial monomial ideal of the Lawrence ideal $J_L$ is
  a radical ideal.
\item{\rm (c)} Every minimal generator of $J_L$ is a difference of two
  squarefree monomials.
\item{\rm (d)} The lattice $L$ is the image of an integer matrix $B$ with
  linearly independent columns, such that all maximal minors of $B$
  lie in the set $\set{ 0, 1, -1}$.
\item{\rm (e)} The lattice $L$ is the kernel of an integer matrix $A$ with
  linearly independent rows, such that all maximal minors of $A$ lie
  in $\set{ 0, m, -m}$ for some integer $m$.
\item{\rm (f)} The quotient ring $S/J_L$ is a normal domain.

\ref{unimod} is proved in \ref{unimodularity}.
If any (and thus all) of these
six equivalent conditions for $L$ holds, then we say that the
Lawrence ideal $J_L$ is {\it unimodular\/}. A first example 
is the ideal of $2 \times 2$-minors of a
$2 \times n$-matrix of indeterminates:
$$
J_L \quad = \quad I_2 \pmatrix{ x_1 & x_2 & x_3 & \cdots & x_n \cr
  y_1 & y_2 & y_3 & \cdots & y_n \cr} \eqno (1.1) $$
Here $L$ is the kernel of 
$A = (\, 1 \, \,1 \, \,1 \, \, \cdots \, 1 \,)$, 
or the image of the matrix $B$ whose rows are $\ee_i -
\ee_{i+1}$, $i\in\set{1,\ldots,n-1}$, the differences of consecutive
unit vectors in $\RR^n$.  The minimal resolution of (1.1) is an
Eagon-Northcott complex, whose polyhedral model is the hypersimplicial
complex of Gel'fand and MacPherson [BS, Ex.~3.15].

In \ref{unimodularity} we introduce an
infinite periodic hyperplane
arrangement $\cHH_L$ whose vertices are the elements of $L$.
It is shown in \ref{main}  that this arrangement
supports the minimal free resolution of $J_L$, and coincides with the hull
complex of $J_L$. In \ref{mono} we prove that this resolution is universal
in the sense that it is stable under
all Gr\"obner deformations. In particular, 
all initial ideals of $J_L$ have
the same Betti numbers as $J_L$, and their minimal resolutions
are also cellular. We also construct cellular resolutions for
the monomial ideals defined by the fibers of $L$.
In \ref{examples} we discuss Lawrence ideals associated
with directed graphs and we present open combinatorial problems.
In \ref{diagonal} we reinterpret
Lawrence ideals in terms of the Audin-Cox homogeneous coordinate ring,
and we generalize Beilinson's spectral sequence 
from projective space to unimodular toric varieties.

The authors would like to thank the Mathematical
Sciences Research Institute in Berkeley for its support while part of
this paper was being written. All three authors were 
partially supported by the NSF during the preparation of this work.

\section{unimodularity} {Unimodularity and an infinite hyperplane arrangement}

\noindent
In this section we establish some basic facts about unimodular
lattices and their Lawrence ideals. We start out by proving
the equivalences stated in \ref{introd}.

\vskip .1cm

\noindent {\bf Proof of \ref{unimod}.}
A monomial ideal is radical if and only if it is
generated by squarefree monomials.  Clearly, the first term of a
binomial $\,\xx^\aa \yy^\bb - \xx^\bb \yy^\aa\,$ is squarefree if and
only if the second term is squarefree. The equivalence of (a), (b)
and (c) follows directly from \ref{lawrencegens}.

In an exact sequence of free abelian groups,
$$
  0 \,\, \longrightarrow \,\,\ZZ^{n-d}
\,\,\, {\buildrel B \over \longrightarrow }  \,\,\, \ZZ^n
\,\,\, {\buildrel A \over \longrightarrow }  \,\,\, \ZZ^d,
\eqno (2.1) $$
each $(n \! - \! d) \times (n \! - \!d)$-minor
of $B$  is equal to the complementary
$d \times d$-minor of $A$, up to a global constant $m$. Since
the cokernel of $B$ is torsion free,
the maximal minors of $B$ are integers
with no common factor. This implies $m \in \ZZ$
and the equivalence of  (d) and (e), for
$\, L = \ker(A) = \im(B)$.

The condition (e) is precisely the defining condition given in [Stu,
\S 8, page 70] for a matrix $A$ to be unimodular. A matrix $A$ is
unimodular if and only if its Lawrence lifting $\Lambda(A)$ is
unimodular. Following [Stu, \S 7, page 55], the {\it Lawrence lifting\/} of
the $d \times n$-matrix $A$ is obtained by appending the zero $d \times
n$-matrix ${\mathbf 0}_{d,n}$ and two copies of the $n \times n$-identity
matrix ${\mathbf I}_n$ as follows: $$ \Lambda(A) \quad = \quad
\pmatrix{ A\hfill &   {\mathbf 0}_{d,n}\hfill \cr
          {\mathbf I}_n\hfill & \, {\mathbf I}_n \hfill\cr}. $$
Hence the equivalence of (b) and (e) is a reformulation
of [Stu, Remark 8.10].

The conjunction of (b) and (e) implies property (f), namely, the
lattice $L$ is the kernel of an integer matrix if and only if $J_L$ is
a prime ideal, and (b) implies normality by [Stu, Proposition 13.15].
To complete the proof, it suffices to show that (f) implies (c).
Suppose that (c) is false, i.e., the ideal $J_L$ has a minimal
generator $\,\xx^\aa \yy^\bb - \xx^\bb \yy^\aa\,$ whose terms are not
squarefree. We may assume that this generator is a circuit, which
means that it has minimal support. By setting all pairs of variables not
appearing in this circuit to zero, we reduce to the case where $J_L$
is a principal ideal. But an affine binomial hypersurface is normal if
and only if at least one of its monomials is squarefree. This
completes our proof of \ref{unimod}.  \Box

Here is another, more invariant, formulation of the unimodularity condition.

\proposition{coordproj} A lattice $L$ is unimodular if and only if,
for every projection $\pi(L)$ of $L$ to a coordinate sublattice $\ZZ^r
\subset \ZZ^n$, the quotient $\ZZ^r/\pi(L)$ is torsion free.

\proof We will use the notation of the proof of \ref{unimod}.
We first prove the ``if'' direction.  Set $m = n-d$ and consider
any nonsingular $m \times m$-submatrix $C$ of $B$. The image of $C$ in $\ZZ^m$
equals $\pi(L)$ for the corresponding coordinate projection $\pi : \ZZ^n
\rightarrow \ZZ^m $.  The finite abelian group $\ZZ^m/\pi(L)$ is torsion free
if and only if it is zero.  Hence $\pi(L) = \ZZ^m$ and we conclude that the
determinant of the $m \times m$-submatrix of $B$ under consideration is
either $1 $ or $-1$.

For the ``only-if'' direction, suppose that $L$ is unimodular.  
Every coordinate projection $\pi(L)$ of $L$ is unimodular as
well. This can be seen by choosing  $A$ (resp.~$B$) to have
appropriate unit vectors among its columns (resp.~rows).  Thus it
suffices to show that the group $\ZZ^n/L$ is torsion free.   But this
follows from the exact sequence (2.1), which identifies $\ZZ^n/L$ with
the image of the matrix $A$, showing
 $\ZZ^n/L \simeq \ZZ^d$.  \Box

Let $L$ be any $m$-dimensional sublattice of $\ZZ^n$, and write $\RR
L$ for the linear subspace of $\RR^n$ spanned by $L$.  We denote by
$\cHH_L$ the {\it affine hyperplane arrangement\/} in $\RR L$ obtained by
intersecting $\RR L$ with all lattice translates of the coordinate
hyperplanes in $\RR^n$. These are the hyperplanes $\set{x_i = j}$ for
$\, 1\le i \le n $ and $ j \in \ZZ $.  Thus $\cHH_L$ is an infinite
$m$-dimensional hyperplane arrangement in the vector space $\RR L$.

It is convenient to embed $\cHH_L$ as a hyperplane arrangement in the
Euclidean space $\RR^m$ as follows. Let $B$ be an integer $n \times
m$-matrix such that  $L=\im(B)$ as in part (d) of \ref{unimod}.  Write
$\bb_i \in \ZZ^m$ for the $i$-th row vector of $B$.  Then
$\cHH_L$ is isomorphic to the infinite arrangement in $\RR^m$ consisting 
of the hyperplanes
$$\, H_{ij} \,=\, \setdef{\xx \in \RR^m}{\bb_i \cdot \xx = j} \quad
\hbox{ for all $i \in\set{1,2,\ldots,n}$ and all $j \in \ZZ$} .$$

\proposition{nonewvertices} Each lattice point in $L$ is a vertex of
the affine hyperplane arrangement $\cHH_L$. There are no
additional vertices in $\cHH_L$ if and only if $L$ is unimodular.

\proof We identify $L$ with $\ZZ^m$ via the matrix $B$ and thus consider
$\cHH_L$ as the hyperplane arrangement $\cHH_L=\set{H_{ij}}$ in $\RR^m$.
The intersection point of $m$ linearly independent such hyperplanes,
$$
H_{i_1 j_1} \, \cap \, H_{i_2 j_2} \, \cap \, \cdots \,
\cap H_{i_m  j_m} \quad = \quad \{ \xx \}, $$
is defined by the linear equations
$\, \bb_{i_\nu} \cdot \xx \, = \,j_\nu \,$ for $\nu = 1,\ldots, m$.
The point $\xx\in\RR^m$ has integer coordinates for all
$j_1,\ldots,j_m \in \ZZ$ if and only if
$\det(\bb_{i_1},\ldots,\bb_{i_m}) = \pm 1$.  By part (d) in
\ref{unimod}, this means that the lattice $L$ is unimodular.  The
first assertion holds because each $\xx \in \ZZ^m$ can be expressed
as such an intersection.  \Box

>From now on we assume that  $L$  is a unimodular
sublattice of $\ZZ^n$.
Recall from [Stu, Proposition 8.11] that the
Graver basis of the Lawrence ideal $J_L$ is given precisely by the
circuits of the lattice $L$. (The {\it circuits\/} of $L$ are the
primitive vectors in $L$ whose supports are minimal 
with respect to inclusion.)
We view $\cHH_L$  as an infinite regular cell complex,
equipped with an action by the abelian group $L$.

\lemma{circuit} Two vertices $\aa$ and $\bb$ of the arrangement $\cHH_L$
are connected by an edge if and only if their difference $\aa - \bb$
is a circuit of the unimodular lattice $L$.

\proof Since $L$ acts transitively on the vertices of $\cHH_L$ we may
assume that $\bb = \zero$. Our assertion states that $\aa$ is a
circuit if and only if $\{\zero,\aa\}$ forms an edge in $\cHH_L$.
This holds because the circuits in the subspace $\RR L$ of $\RR^n$ are
computed by the rule
$$  H_{i_1 0} \, \cap \,
 H_{i_2 0} \, \cap \, \cdots \, \cap
 H_{i_{m-1} 0} \quad = \quad \RR \aa  $$
for all possible increasing sequences of
indices $\,1 \leq i_1 < i_2 < \cdots < i_{m-1} \leq n $.
\Box

In this section we have introduced a family of hyperplane arrangements
$\cHH_L$ whose vertices form a lattice $L$.  Such arrangements appear
in many parts of the mathematical literature; for example, see
[BLSWZ]. The group $L$ acts on the faces
of the arrangement $\cHH_L$ with finitely many orbits, and we shall be
interested in the quotient complex $\cHH_L/L$.
Nontrivial examples will be presented in
\ref{examples}. First we reinterpret $\cHH_L$ and $\cHH_L/L$ as
minimal free resolutions in the sense of commutative algebra. This 
will be done in the next two sections by labeling the faces of
$\cHH_L$ with Laurent monomials, following the general recipes in
[BS].

Our warmup example (1.1) is the case where $L$ is spanned by the
vectors $\ee_i - \ee_{i+1}$, $i\in\set{1,\ldots,n-1}$, the differences
of the consecutive unit vectors in $\RR^n$. Thus $J_L$ is here the ideal
of $2 \times 2$-minors, and $\cHH_L$ is the {\it hypersimplicial
arrangement\/}. The quotient complex ${\cal H}_L/L$ has $n-1$ distinct
maximal faces, called {\it hypersimplices\/}.
The hypersimplicial arrangement for $n=3$ is depicted in the following
figure:

\capdraw{100}{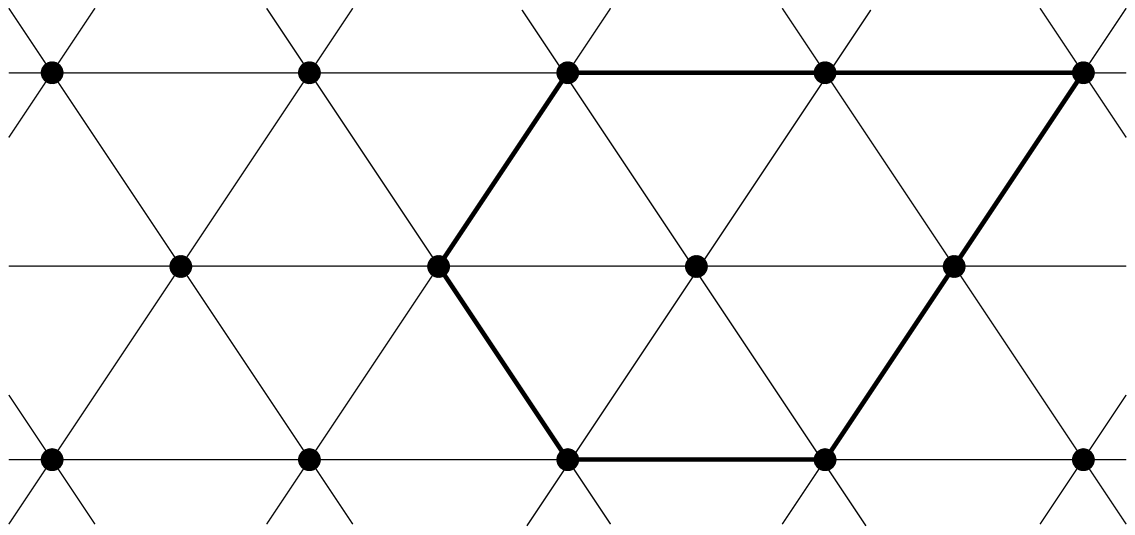}{
\swtext{.075}{.875}{\swpad{2pt}{2pt}{${{y_1^2y_2^{}x_3^3}
  \over{x_1^2x_2^{}y_3^3}}$}}
\swtext{.305}{.875}{\swpad{2pt}{2pt}{${{y_1^{}y_2^{}x_3^2}
  \over{x_1^{}x_2^{}y_3^2}}$}}
\swtext{.530}{.875}{\swpad{2pt}{2pt}{${{y_2^{}x_3^{}}\over{x_2^{}y_3^{}}}$}}
\swtext{.760}{.875}{\swpad{2pt}{2pt}{${{x_1^{}y_2^{}}\over{y_1^{}x_2^{}}}$}}
\swtext{.990}{.875}{\swpad{2pt}{2pt}{${{x_1^2y_2^{}y_3^{}}
 \over{y_1^2x_2^{}x_3^{}}}$}}
\swtext{.190}{.505}{\swpad{2pt}{2pt}{${{y_1^2x_3^2}\over{x_1^2y_3^2}}$}}
\swtext{.415}{.505}{\swpad{2pt}{2pt}{${{y_1^{}x_3^{}}\over{x_1^{}y_3^{}}}$}}
\swtext{.635}{.505}{\swpad{4pt}{4pt}{$1$}}
\swtext{.870}{.505}{\swpad{2pt}{2pt}{${{x_1^{}y_3^{}}\over{y_1^{}x_3^{}}}$}}
\swtext{.075}{.139}{\swpad{2pt}{2pt}{${{y_1^3x_2^{}x_3^2}
 \over{x_1^3y_2^{}y_3^2}}$}}
\swtext{.305}{.139}{\swpad{2pt}{2pt}{${{y_1^2x_2^{}x_3^{}}
  \over{x_1^2y_2^{}y_3^{}}}$}}
\swtext{.530}{.140}{\swpad{2pt}{3pt}{${{y_1^{}x_2^{}}\over{x_1^{}y_2^{}}}$}}
\swtext{.760}{.139}{\swpad{2pt}{2pt}{${{x_2^{}y_3^{}}\over{y_2^{}x_3^{}}}$}}
\swtext{.990}{.139}{\swpad{2pt}{2pt}{${{x_1^{}x_2^{}y_3^2}
  \over{y_1^{}y_2^{}x_3^2}}$}}
 }
{The hypersimplicial arrangement for $n=3$.}

Its hypersimplices are up-triangles and down-triangles.
Each vertex is labeled by a Laurent monomial.
The finite quotient complex  $\cHH_L/L$ is  a torus, subdivided
by one vertex, three edges and two $2$-cells,
one  up-triangle and  one down-triangle.

\section{main} {From hyperplane arrangements to minimal free resolutions}

\noindent
We fix a  unimodular sublattice $L$ of dimension $m$ in $\ZZ^n$.
The Laurent polynomial ring
$T=k[x_1^{\pm 1},\ldots,x_n^{\pm 1}, \, y_1^{\pm 1},\ldots,y_n^{\pm1}]$
is a module over the polynomial ring
$S =  k[x_1,\ldots,x_n, \, y_1,\ldots,y_n]$. We consider
the following monomial $S$-submodule of $T$:
$$ M_L \quad := \quad S \cdot \setdef{\xx^\aa \yy^{-\aa}}{\aa \in L}
\quad \subset \quad T.$$
Each lattice point $\aa = (a_1,\ldots,a_n)$ in $L$
is a vertex of the arrangement $\cHH_L$, and we label
that vertex of $\cHH_L$ with the corresponding generator of $M_L$, namely,
$$  \xx^\aa \yy^{-\aa}  \quad := \quad
x_1^{a_1} x_2^{a_2}  \cdots x_n^{a_n}
y_1^{-a_1} y_2^{-a_2}  \cdots y_n^{-a_n}. $$
Each face $F$ of $\cHH_L$ is then labeled by the least common multiple
$m_F$ of the labels of its vertices. (The {\it least common multiple\/}
of a set of Laurent monomials is the Laurent monomial
whose exponents are the coordinatewise maxima of the given exponents.)
The labeled cell complex $\cHH_L$ defines a
complex of free $\ZZ^{2n}$-graded $S$-modules
$$\FF_{\cHH_L} \,\, = \,\,\bigoplus_{F\in\cHH_L \atop F\ne\emptyset} S(
-m_F ),$$
where the summand $S(-m_F)$ has homological degree $\dim(F)$.
The differential of the complex $\FF_{\cHH_L}$ is the homogenized differential
of the cell complex $\cHH_L$, defined by
$$\partial(F) \quad = \,\sum_{F'\subset F \atop \codim(F',F)=1}
\epsilon(F,F') \cdot {m_F \over m_{F'}}\cdot F',\qquad
\hbox{for faces} \quad  F \,\,\,\hbox{of} \,\,\, \cHH_L.$$
Here  $\epsilon(F,F')$  is either $+1$ or $-1$, indicating
the orientation of $F'$ in the boundary of $F$.
See [BS, \S1] for details on this construction.
The complex $(\FF_{\cHH_L},\partial)$ is
not $S$-finite, but has finite length $m=\rank(L)$.

As an example consider the complex of free
$k[x_1,x_2,x_3,y_1,y_2,y_3]$-modules defined by the hypersimplicial
arrangement in \ref{fig:hyper}. The edge $E$ connecting the module
generators $\,m_1 = {x_2 y_3 \over y_2 x_3}\,$ and $\,m_2 = {x_1 y_3
  \over y_1 x_3}\,$ is labeled by their least common multiple, which
is the Laurent monomial $\,m_E = {x_1 x_2 y_3 \over x_3}$.  This edge
$E$ represents the first syzygy $( x_1 y_2) \cdot m_1 - (y_1 x_2)
\cdot m_2$ of the monomial submodule $M_L$.  

\theorem{mainthm}  The complex $(\FF_{\cHH_L},\partial)$
is a minimal $\ZZ^{2n}$-graded free $S$-resolution of
the lattice module $M_L$.

\proof
The complex $(\FF_{\cHH_L},\partial)$ consists of free $S$-modules and
is clearly $\ZZ^{2n}$-graded. To show that it is a resolution, we
apply the exactness criterion in [BS, Proposition 1.2] to the labeled
cell complex $X = \cHH_L$.  For any $(\aa,\bb) \in \ZZ^{2n}$ we
consider the subcomplex $\,X_{\leq (\aa,\bb)} \,$ consisting of all
faces $F$ of the arrangement $X$ whose label $\,m_F = \xx^\cc \yy^\dd
\,$ satisfies the coordinatewise inequalities $\,\cc \leq \aa \,$ and
$\,\dd \leq \bb $.  We shall prove that $X_{\leq (\aa,\bb)}$
is contractible, by identifying this subcomplex with a convex polytope
in $\RR L$.  For instance, the marked pentagon in
\ref{fig:hyper} is  $\,X_{\leq (2,1,1,\, 1,1,1)}$.

Our labeled hyperplane arrangement can be described as follows.  For
$\uu \in \RR^n$ we write $\,(\ceil{\uu}, \floor{\uu})\,$ for the vector in
$\ZZ^{2n}$ obtained by rounding up and down each coordinate of $\uu$.  For instance,
if $\uu = ( -2/5,-1, 7/5)\,$ then $\,(\ceil{\uu}, \floor{\uu}) = (0,-1,2 , -1,
-1,1) $. Two vectors $\uu$ and $\vv$ in the subspace $\RR L$ of
$\,\RR^n \,$ are called ${\it equivalent\/}$ if $\,(\ceil{\uu}, \floor{\uu}) \, =
\,(\ceil{\vv}, \floor{\vv}) $.  The resulting equivalence classes on $\RR L$ are
the (relatively open) faces of $X = {\cal H}_L$. The face $F$ containing $\uu \in
\RR L$ is labeled by the vector $\, \,(\ceil{\uu}, -\floor{\uu})$, or by the
corresponding Laurent monomial
$m_F=\xx^{\ceil{\uu}} \yy^{-\floor{\uu}}$.

For $\aa,\bb \in \ZZ^n$ we consider the following subset of
the subspace $\,\RR L \,$ in $\, \RR^n$:
$$
\setdef{\uu \in \RR L}{
\ceil{\uu}  \leq \aa \,\,\, \hbox{and} \,\,\,
-\floor{\uu} \leq \bb}
\quad = \quad
\setdef{ \uu \in \RR L}{
 - \bb \leq \uu \leq \aa}. \eqno (3.1{\rm a})
$$
This is a convex polytope with facet hyperplanes taken from $\cHH_L$.
By the construction in the previous paragraph, the complex $\,X_{\leq
  (\aa,\bb)} \,$ is a polyhedral subdivision of the convex polytope
identified in (3.1a).  Therefore $\,X_{\leq (\aa,\bb)} \,$ is
contractible, and, by [BS, Proposition 1.2], the complex
$(\FF_{\cHH_L},\partial)$ is exact over $S$.

It follows from the description of the labeled complex $
X = \cHH_L$ in the second-to-last paragraph that distinct faces
$F$ and $F'$ of $X$ have distinct labels $m_F \not= m_{F'}$. This
shows that the resolution $(\FF_{\cHH_L},\partial)$ is  minimal
(compare [BS, Remark 1.4]).  \Box

\corollary{ZeroOne}
$\!\! \!\! $ The $\ZZ^{2n}$-graded Betti numbers of 
$\,M_L\,$ are $0$ or $1$.
  \Box


We retain the following identification for the rest of the paper:
$$
X_{\leq (\aa,\bb)} \quad = \quad
\setdef{\uu \in \RR L}{- \bb \leq \uu \leq \aa}
\qquad \hbox{for} \quad \aa ,\bb \in \ZZ^n. \eqno (3.1{\rm b})
$$
In the next section we shall further identify these polytopes with the
{\it fibers\/} of the Lawrence ideal $J_L$, i.e., with the congruence
classes modulo $J_L$ of monomials in $S$.

Note that the polytope $\,
X_{\leq (\aa,\bb)} \,$ is itself a (closed) face of the arrangement $X
= \cHH_L$ if and only if there exists a vector $\uu \in \RR L $
such that $\,(\aa,\bb) = (\ceil{\uu}, -\floor{\uu})$.
Suppose that this holds. Then each coordinate of $\,\aa + \bb \,$ is
either $0$ or $1$, and, taking any vertex of $\,X_{\leq (\aa,\bb)}$,
we get a lattice point $\vv \in L$ with $\,-\bb \leq \vv \leq \aa $. By
translating our face with the lattice vector $\vv$, we obtain now the
following conclusion:

\proposition{quotientL} Modulo the action by the unimodular lattice
$L$, each face of ${\cal H}_L$ has the form $\,X_{\leq (\aa,\bb)}$,
where $\aa \in \set{0,1}^n$ and $\bb\in \set{0,1}^n$ have disjoint
support.  \Box

We next consider the quotient complex $\cHH_L/L$, which is formally
defined as the face poset of $\cHH_L$ modulo the action
by the lattice $L$. \ref{quotientL} implies:

\corollary{isFinite}
The number of faces of the quotient complex $\cHH_L/L$ is finite.
  \Box

We apply the algebraic quotient construction in [BS, Section 3] to
this quotient complex. The labeling of the $\ZZ^{2n}$-graded cell
complex $X = {\cal H}_L$ is consistent with the action by the
following lattice which is canonically isomorphic to $L$,
$$ \Lambda(L) \quad := \quad
\setdef{ (\uu,-\uu) \in \ZZ^{2n}}{ \uu \in L}.
\eqno (3.2) $$
Following [BS, Lemma 3.5], the complex $(\FF_{\cHH_L},\partial)$ has
the structure of a complex of $\ZZ^{2n}$-graded free modules over the
group algebra $\, S[L]$.  The rank of $\FF_{\cHH_L}$ over $S[L]$
equals the number of faces of $\cHH/L$, which is finite by
\ref{isFinite}. Now the functor in [BS, Theorem 3.2]
defines an equivalence of categories between the category of
$\ZZ^{2n}$-graded $S[L]$-modules and the category of
$\ZZ^{2n}/\Lambda(L)$-graded $S$-modules. Applying this functor to the
$S[L]$-complex $(\FF_{\cHH_L},\partial)$ we obtain the cellular
quotient complex $\,(\FF_{\cHH_L/L},\partial)\,$ which is a complex of
$\ZZ^{2n}/\Lambda(L)$-graded free $S$-modules.  From \ref{mainthm} and
[BS, Corollary 3.7] we conclude that $\,(\FF_{\cHH_L/L},\partial)\,$
is the minimal free resolution of $S/J_L$ over $S$. This concludes the proof
of the following theorem.

\theorem{unimod-toric-syz}
The quotient complex $\cHH_L/L$ of the hyperplane
arrangement modulo $L$ supports the minimal $S$-free
resolution of the unimodular Lawrence ideal $J_L$.  \Box

\corollary{ZeroOneagain} The minimal free resolution of
the unimodular Lawrence ideal $J_L$ is independent of
the characteristic of the base field $k$. The number of minimal
$i$-th syzygies of $S/J_L$ equals the number of $i$-dimensional
faces of the quotient complex $\cHH_L/L$, and the Betti numbers
of $S/J_L$  in the $\ZZ^{2n}/\Lambda(L)$-grading
are all $0$ or $1$.
  \Box


To write down the matrices in the minimal cellular resolution
$\,(\FF_{\cHH_L/L},\partial)\,$ of a unimodular Lawrence ideal $J_L$,
one must select a fundamental domain of $\cHH_L$ modulo $L$ and identify
the cover relation in the poset of faces of $\cHH_L/L$.
In higher dimensions it is convenient to use the monomial ideals
in \ref{mono} for that purpose, but in
dimensions $2$ and $3$ we can do the identifications directly on the
picture. For instance, if
$\,L = \ker(1 \,\, 1 \,\, 1 )\,$
and thus $J_L$ is the ideal of $2 \times 2$-minors
of a generic $2 \times 3$-matrix, then the resolution
$\,(\FF_{\cHH_L/L},\partial)\,$ is supported by the hypersimplicial
arrangement in $\RR^2$ from \ref{fig:hyper}. 
Here the quotient complex
$\cHH_L/L$ consists of one vertex, three edges and two triangles,
which are glued to form a torus. Hence $S/J_L$ has
one generator, three first syzygies and two second 
syzygies:

$$
0
\longrightarrow  S^2
\rightarrowmat{2pt}{4pt}{x_1 & y_1 \cr x_2 & y_2 \cr x_3 & y_3 \cr} S^3
\rightarrowmat{2pt}{4pt}{x_2 y_3 \! - \! x_3 y_2 & x_3 y_1 \! - \! x_1 y_3 &
x_1 y_2 \! - \! x_2 y_1 \cr} S
\longrightarrow S/J_L\longrightarrow 0
$$

\vskip .1cm

We now show that the minimal cellular resolution of $S/J_L$
equals the hull resolution introduced in [BS, Section 3].
For  $\aa$, $\bb\in\ZZ^n$ and $t > 0$, 
we write
$$(t^\aa,t^\bb) \quad = \quad
(t^{a_1},\ldots,t^{a_n},t^{b_1},\ldots,t^{b_n}) \, \in \,
\RR^{2n}$$
Recall that $\hull(M_L)$ is the complex of
bounded faces of the polyhedron $P_t$ for large $t$, where $P_t$ is the 
convex hull of
$\setdef{(t^\aa,t^\bb)}{\xx^\aa \yy^\bb \in M_L}\subset \RR^{2n}$. The 
vertices of $P_t$ correspond to the minimal generators of $M_L$, and 
the {\it hull resolution\/} of
$S/J_L$ is the cellular resolution supported on
$\hull(M_L)/L$. We shall use the following lemma.

\lemma{hardanalysis} Let $\set{a, b}$ and  $\set{c,d}$ be pairs of integers so
$a+b=c+d$, and let $t>0$, $t\ne 1$. If $\abs{a-b} = \abs{c-d}$, then $t^a + t^b =
t^c + t^d$. If $\abs{a-b} > \abs{c-d}$, then $t^a + t^b > t^c + t^d$.

\proof
If $\abs{a-b} = \abs{c-d}$, then $\set{a,b} = \set{c,d}$ as multisets.
If $\abs{a-b} > \abs{c-d}$, suppose to be definite that $a > c \ge d > b$. 
In view of $a-c=d-b$, we have
$\,
t^a-t^c-t^d+t^b
\;=\;
t^c(t^{a-c}-1)-t^b(t^{d-b}-1)
\;=\;
t^b(t^{c-b}-1)(t^{a-c}-1)
\;>\;
0 $.
\Box

\theorem{hullthm}
The hull resolution of $S/J_L$ agrees with the minimal free resolution.

\proof
We show that
$\hull(M_L)$ and $\cHH_L$ agree as cell complexes, using the fact that $\hull(M_L)$
consists of those faces of $P_t$ supported by strictly positive inner normals.

Let $F$ be a face of $\cHH_L$. Then $F = X_{\leq (\ceil{\uu}, -\floor{\uu})}$ for
some $\uu \in \RR L$.
In other words, the vertices of $F$ are precisely those elements $\aa\in L$ so
$\floor{u_i} \le a_i \le \ceil{u_i}$ for each $i$.
Let $\vv=(t^{-\ceil{\uu}}, t^{\floor{\uu}})$. If $\aa$ is a vertex
of $F$, then $\vv\cdot(t^\aa,t^{-\aa}) = n + m +(n-m)t^{-1}$, where $m$ is the
number of coordinates in which $\uu$ is an integer. Suppose that $\bb\in L$ is not a
vertex of $F$, and let $t>1$. Write $\vv_i = (t^{-\ceil{u_i}}, t^{\floor{u_i}})$, so
$\vv\cdot(t^\aa,t^{-\aa}) = \sum_{i=1}^n \vv_i \cdot (t^{a_i},t^{-a_i})$.
For each $i$ we have
$\vv_i\cdot (t^{b_i},t^{-b_i}) \ge
\vv_i\cdot (t^{a_i},t^{-a_i})$,
with equality if and only if $\floor{u_i} \le b_i \le \ceil{u_i}$.
This follows by applying \ref{hardanalysis} to the pairs 
$\set{b_i-\ceil{u_i}, \floor{u_i}-b_i}$ and
$\set{a_i-\ceil{u_i}, \floor{u_i}-a_i}$. 
Thus $\vv\cdot(t^\bb,t^{-\bb}) > \vv\cdot(t^\aa,t^{-\aa})$, so $\vv$ supports $F$ as
a face of $\hull(M_L)$.

Let $\aa$ and $\bb$ be two vertices of $\cHH_L$ which do not belong to a common
face of $\cHH_L$, and let $F = X_{\leq (\ceil{\uu},
-\floor{\uu})}$ be the face of $\cHH_L$ determined by $\uu=(\aa+\bb)/2$. Then
$\abs{a_i-b_i} \ge 2$ for some $i$, for otherwise $\aa$ and
$\bb$ would both belong to $F$ by \ref{quotientL}. Let $\cc$ be any vertex of $F$,
and let
$\dd=\aa+\bb-\cc$. Note that for each $j$, $\set{c_j,
d_j}=\set{\floor{u_j}, \ceil{u_j}}$ as multisets, so $\abs{c_j-d_j}\le 1$.
Then
$\aa+\bb=\cc+\dd$, and $\abs{a_j-b_j} \ge \abs{c_j-d_j}$ for each $j$, with strict
inequality for $j=i$. Let
$t>1$. Applying
\ref{hardanalysis} to the pairs $\set{a_j,b_j}$ and $\set{c_j,d_j}$,
$$t^{a_j}+t^{b_j} \ge t^{c_j}+t^{d_j} \quad\mtext{and}\quad
t^{-a_j}+t^{-b_j} \ge t^{-c_j}+t^{-d_j}$$
for each $j$, with strict inequalities for $j=i$. 
Let $\pp$ be the midpoint of the line segment from $(t^\aa,t^{-\aa})$ to
$(t^\bb,t^{-\bb})$, and
let $\qq$ be the midpoint of the line segment from $(t^\cc,t^{-\cc})$ to
$(t^\dd,t^{-\dd})$. We have shown that $\pp-\qq$ is a nonzero, nonnegative
vector. Therefore, the point $\pp$ cannot lie on any face of
$\hull(M_L)$, because $\vv\cdot\pp > \vv\cdot\qq$ for any strictly positive vector
$\vv$. Thus $\aa$ and
$\bb$ cannot belong to a common face of $\hull(M_L)$.
We conclude that the cell complexes
 $\hull(M_L)$ and $\cHH_L$ are equal.
\Box

\section{mono} {Fiber monomial ideals and initial monomial ideals}

\noindent With any lattice ideal in a polynomial ring one can
associate two families of monomial ideals.
First, there are the {\it initial monomial ideals\/},
with respect to various term orders. Their minimal free resolutions
can always be lifted, by Schreyer's construction
in Gr\"obner basis theory, to a (possibly nonminimal)
resolution of the lattice ideal. Second, we have
the {\it fiber monomial ideals\/}, which are generated by the
fibers, that is, the equivalence classes of monomials
modulo the lattice ideal [PS, Section 2]. It follows from
the results in [BS, Section 3] that the resolution of the lattice ideal
is always determined by the resolution of a large enough fiber ideal.

In this section we make these results precise
for the unimodular Lawrence case.
Fix a unimodular sublattice $L\subset\ZZ^n$. Two monomials $m$ and $m'$
in $S = k[x_1,\ldots,x_n,y_1,\ldots,y_n]$ are considered
{\it equivalent\/} if $m - m'$ lies in the unimodular
Lawrence ideal $J_L$. The equivalence classes are finite
and called the {\it fibers\/} of $J_L$.
For given $\aa,\bb \in \NN^n$, let
$\fib(\aa,\bb)$ denote the fiber
of the monomial $\,\xx^\aa \yy^\bb$.
We shall identify the fibers with the
lattice points in the polytopes of the form (3.1b):

\lemma{lawrence-fibers}
Let $X = \cHH_L$ be the labeled cell complex introduced in
\ref{main}. Then the map 
$\,\,\, \phi \, : \, {\fib}(\aa,\bb) \rightarrow
X_{\leq (\aa,\bb)} \, \cap \, L \, , \,\,
\xx^\cc \yy^\dd \mapsto \aa - \cc  \,\,\,$
is a bijection.

\proof
A monomial $\,\xx^\cc \yy^\dd \,$ in $S$ lies in ${\fib}(\aa,\bb)$ if
and only if the vector $\,(\aa - \cc, \bb - \dd)\in\ZZ^{2n}$ lies
in the lattice $\Lambda(L)$ defined in (3.2).  The latter condition
means that $\, \uu = \aa - \cc \,$ lies in $L$ and $\, \bb + \uu = \dd
$.  Thus $\,\uu \,$ is a vertex of $X = \cHH_L$ and its label
$$ \xx^\uu \yy^{-\uu} \quad = \quad
\xx^{\aa - \cc} \yy^{\bb - \dd} \quad  \hbox{divides} \quad
\xx^{\aa} \yy^{\bb } , $$
which shows that $\uu$ is actually a vertex of
$X_{\leq (\aa,\bb)}$. Hence the map $\phi$
is well-defined. To see that it is a bijection
we note that the inverse map is given by
$\,\phi^{-1} (\uu) \,=\, \xx^{ \aa - \uu} \yy^{ \bb + \uu} $.
This map is an analogue of the bijection in
[PS, (2.1)]. \Box

Let $\, \langle {\fib}(\aa,\bb) \rangle \,$ denote the ideal in $S$
generated by all monomials in the fiber of $\xx^\aa \yy^\bb$. The
vertices of the cell complex $\, X_{\leq (\aa,\bb)} \, $ (considered
as a subcomplex of $\cHH_L$) are labeled with certain Laurent
monomials of the form $\,\xx^{\uu} \yy^{-\uu}$, with $\uu \in
L$. Consider their preimages under the bijection $\phi$, and let $\,
Y_{(\aa,\bb)} \, $ denote the same cell complex as $\, X_{\leq
(\aa,\bb)} \, $ but with the vertices labeled by the monomials in
${\fib}(\aa,\bb)$.  Thus the vertex with label $\,\xx^{\uu}
\yy^{-\uu}\,$ in $ \,X_{\leq (\aa,\bb)} \, $ becomes the vertex with
label $\,\xx^{\aa-\uu} \yy^{\bb+\uu} \,$ in $\,Y_{(\aa,\bb)}$.  The
labeled cell complex $\,Y_{(\aa,\bb)}\,$ gives rise to a
$\ZZ^{2n}$-graded complex of free $S$-modules as in [BS, Section 1]. 
This complex is always exact and minimal:

\theorem{fiberRes}
Let $L$ be a unimodular sublattice of $\ZZ^n$ and 
$\aa,\bb \in \NN^n$.  The labeled cell complex  $Y_{(\aa,\bb)}$
defines a minimal free resolution of the
monomial ideal $\, \langle  {\fib}(\aa,\bb) \rangle $.

\proof
For any $\cc,\dd \in \NN^n$ we can make the  following
identification of cell complexes
$$ (Y_{(\aa,\bb)})_{\leq (\cc,\dd)} \quad = \quad
X_{\leq ( {\min}\{\bb,\cc-\aa\}, \,
{\min}\{\aa,\dd-\bb \} )}, \eqno (4.1) $$
where ``$\min{}\!$'' refers to the coordinatewise minimum,
and the labels are shifted appropriately.
To see that (4.1) holds, one identifies
both sides with the convex polytope
$$ \setdef{ \uu \in {\RR L}}
{(\zero,\zero) \leq (\aa+ \uu, \bb-\uu) \leq (\cc,\dd)}, $$
together with its natural subdivision.
We conclude that the relevant subcomplexes for the
exactness criterion in [BS, Proposition 1.2] are
all contractible.
As above, in the proof of \ref{mainthm}, distinct faces of
$\,Y_{(\aa,\bb)}\,$ have distinct labels.
Therefore $\,Y_{(\aa,\bb)} \,$
supports a minimal cellular resolution of the fiber ideal
$\,\langle {\fib}(\aa,\bb) \rangle $.
\Box

\example{magicmatrices}
Let $L$ be the corank $1$ lattice of all vectors
in $\ZZ^n$ with coordinate sum zero. As discussed at
the end of \ref{unimodularity}, $X = \cHH_L$ is the
hypersimplicial arrangement in $\RR^{n-1}$, while $J_L$ is the ideal of
$2 \times 2$-minors of the matrix (1.1).
Let $M$ be the ideal generated by the monomials
$\, x_1^{u_1}  x_2^{u_2} \cdots x_n^{u_n}
y_1^{v_1}  y_2^{v_2} \cdots y_n^{v_n}  \,$
where
$$ \pmatrix{
u_1 & u_2 &   \cdots & u_n \cr
v_1 & v_2 &   \cdots & v_n \cr}  \eqno (4.2) $$
runs over all nonnegative integer matrices
having the same row sums and the same column sums.
Then the minimal free resolution of $M$ is cellular
and supported by a suitably relabeled
subcomplex $Y_{(\aa,\bb)}$ of
the hypersimplicial arrangement $X = \cHH_L$.

For instance, let $n=3$ and consider the following monomial ideal
$$ \eqalign{ M \quad = \quad
\langle  \;
    x_2^2 x_3^2 y_1^3 ,  \;
   x_1 x_2 x_3^2 y_1^2 y_2 , \; &
   x_1^2 x_3^2 y_1 y_2^2 , \;
   x_1 x_2^2 x_3 y_1^2 y_3 , \; \cr
   & x_1^2 x_2 x_3 y_1 y_2 y_3 , \;
   x_1^3 x_3 y_2^2 y_3 , \;
   x_1^2 x_2^2 y_1 y_3^2 , \;
   x_1^3 x_2 y_2 y_3^2
\; \rangle                           \cr}
$$
Here the matrices (4.2) have row sums $4,3$ and column
sums $3,2,2$. The minimal free resolution of $M$ is supported on the
complex $Y_{(2,1,1,\, 1,1,1)}$ which is a pentagon:

\capdraw{100}{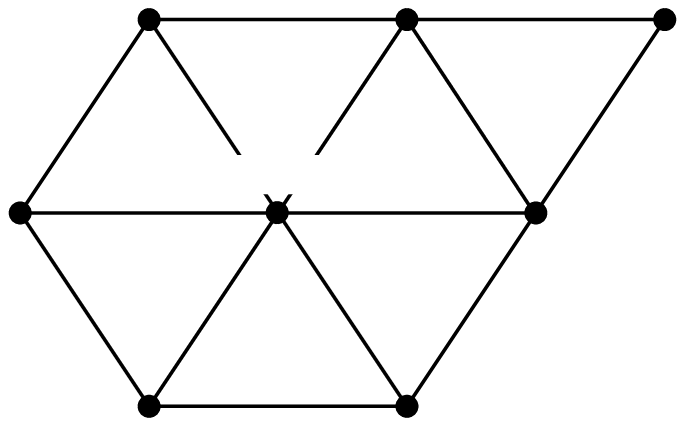}{
\stext{.212}{1}{\spad{3pt}{\small${x_1^2x_2^2y_1^{}y_3^2}$}}
\stext{.595}{1}{\spad{3pt}{\small${x_1^{}x_2^2x_3^{}y_1^2y_3^{}}$}}
\stext{.979}{1}{\spad{3pt}{\small${x_2^2x_3^2y_1^3}$}}
\etext{.000}{.500}{\epad{4pt}{\small${x_1^3x_2^{}y_2^{}y_3^2}$}}
\stext{.403}{.538}{\spad{3pt}{\small${x_1^2x_2^{}x_3^{}y_1^{}y_2^{}y_3^{}}$}}
\wtext{.8080}{.500}{\wpad{4pt}{\small${x_1^{}x_2^{}x_3^2y_1^2y_2^{}}$}}
\etext{.195}{.036}{\epad{4pt}{\small${x_1^3x_3^{}y_2^2y_3^{}}$}}
\wtext{.614}{.036}{\wpad{4pt}{\small${x_1^2x_3^2y_1^{}y_2^2}$}}
}{Minimal free resolution of a fiber in the hypersimplicial arrangement.}

Let $\prec$ be any term order on the polynomial ring $\, S =
k[x_1,\ldots,x_n,y_1,\ldots,y_n]$. Consider the initial monomial ideal
$\,{\ini}_\prec(J_L)\,$ of the unimodular Lawrence ideal $J_L$. We
know from \ref{unimod} that $\ini_\prec(J_L)$ is a squarefree monomial
ideal.  We shall describe a minimal cellular free resolution of
$\,{\ini}_\prec(J_L)\,$ and show that it has the same Betti numbers as
the resolution $\,(\FF_{\cHH_L/L},\partial)\,$ of the Lawrence ideal $J_L$.

Let $\cHH^0_L$ be the set of all faces of the infinite hyperplane
arrangement $\cHH_L$ which contain the origin $\zero \in L$. This is a
finite cell complex which we identify with the central hyperplane
arrangement in $\RR L $ given by the $n$ coordinate hyperplanes $x_i =
0$. The faces of ${\cHH}^0_L$ are cones in $\RR L$ with their apex at
the origin. Under the embedding of $\cHH_L$ in $ \RR^m$ given prior to
\ref{nonewvertices}, the complex ${\cHH}^0_L$ becomes the central
arrangement defined by the hyperplanes $\,H_{10}, H_{20},\ldots,
H_{n0}\,$ in $\RR^m$.  For instance, in the 
example of \ref{fig:hyper}, $\cHH^0_L$ consists of one
$0$-face, six $1$-faces and six $2$-faces.

The term order $\prec$ extends uniquely to a total order on all
Laurent monomials $\xx^{\uu} \yy^{-\uu}$ and hence on all vertices
$\uu$ of $\cHH_L$.  A positive-dimensional cone $F$ in the central
arrangement ${\cHH}^0_L$ is called {\it $\prec$-positive\/} if all
nonzero vertices $\uu \in L $ of the face $F$ satisfy the inequality
$\,1 \prec \xx^{\uu} \yy^{-\uu} $.  We can represent the term order
$\prec$ by a generic hyperplane $H_\prec$ not containing the origin in
$\RR L$, such that the $\prec$-positive cones of ${\cHH}^0_L$ are
precisely those cones which have bounded, nonempty intersection with
$H_\prec$. If $F$ is an $i$-dimensional $\prec$-positive cone in
$\cHH^0_L$ then $ \, F \,\cap \,H_\prec\,$ is an
$(i\! - \! 1)$-dimensional convex polytope.
We write ${\ini}_\prec(\cHH_L)$ for the $(m-1)$-dimensional
labeled cell complex consisting of those (bounded) polytopes
$ \, F \,\cap \,H_\prec$.  Here $ F \,\cap \,H_\prec $ inherits the
label $m_F$ from the face $F$ of $\cHH_L$.  Note that $m_F$ is
a monomial in $S$ since $\,\zero \in F$.

\capdraw{100}{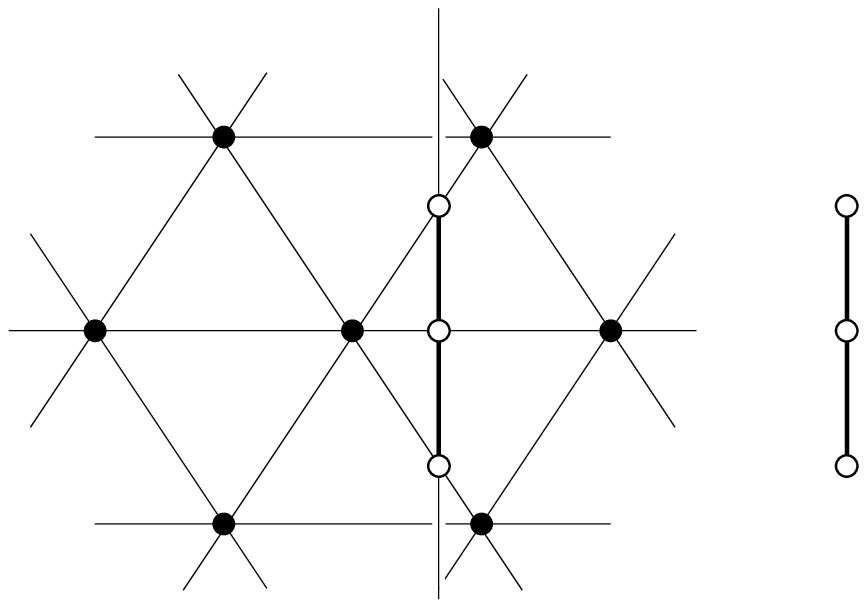}{
\wtext{.507}{.950}{\wpad{3pt}{$H_\prec$}}
\swtext{.290}{.780}{\swpad{2pt}{2pt}{${{y_2x_3}\over{x_2y_3}}$}}
\swtext{.595}{.780}{\swpad{2pt}{2pt}{${{x_1y_2}\over{y_1x_2}}$}}
\swtext{.140}{.460}{\swpad{2pt}{2pt}{${{y_1x_3}\over{x_1y_3}}$}}
\stext{.405}{.480}{\spad{4pt}{$1$}}
\swtext{.745}{.460}{\swpad{2pt}{2pt}{${{x_1y_3}\over{y_1x_3}}$}}
\swtext{.295}{.135}{\swpad{2pt}{2pt}{${{y_1x_2}\over{x_1y_2}}$}}
\swtext{.595}{.135}{\swpad{2pt}{2pt}{${{x_2y_3}\over{y_2x_3}}$}}
\wtext{.995}{.666}{\wpad{3pt}{$x_1y_2$}}
\wtext{.995}{.456}{\wpad{3pt}{$x_1y_3$}}
\wtext{.995}{.230}{\wpad{3pt}{$x_2y_3$}}
\ntext{1.0}{.140}{\npad{3pt}{${\ini}_\prec(\cHH_L)$}}
}{Initial monomial ideal of a unimodular Lawrence ideal.}

\ref{fig:top} shows such a hyperplane $H_\prec$ and the
resulting complex $ \ini_\prec(\cHH_L) $ arising from 
the hypersimplicial complex in \ref{fig:hyper}.
The $1$-dimensional complex $ \ini_\prec(\cHH_L) $ represents the 
minimal $S$-free resolution of the monomial ideal $in_\prec(J_L)$:
$$  0
\longrightarrow
S^2
\rightarrowmat{2pt}{4pt}{x_1 & 0 \cr - x_2 & - y_2 \cr 0 & y_3 \cr} S^3
\rightarrowmat{2pt}{4pt}{ x_2 y_3  &  x_1 y_3  & x_1 y_2  \cr}
S \longrightarrow S/{\ini_\prec(J_L)} \longrightarrow 0 $$

\theorem{initialRes}
The minimal free resolution of the initial monomial ideal
$\,{\ini}_\prec(J_L) \,$  of the unimodular
Lawrence ideal $J_L$
is given by the
labeled cell complex ${\ini}_\prec(\cHH_L)$.

\proof
The $\prec$-positive $1$-dimensional faces $C$ of  ${\cHH}^0_L$
are the circuits of $L$, and their labels
$m_C$ are the minimal generators $\,{\ini}_\prec(J_L) $.
Any higher-dimensional $\prec$-positive face $F$ of ${\cHH}^0_L$
corresponds to a region in the central hyperplane arrangement in $\RR L $ 
given by the coordinates. Its label $m_F$ is a squarefree monomial 
by \ref{quotientL}. We can  identify $m_F$ with the signed vector
in the oriented matroid which represents the region $F$.
Now, every vector in an oriented
matroid is a conformal union of the circuits
[BLSWZ, Proposition 3.7.2]. This implies that
$m_F$ is the least common multiple of the labels $m_C$ of
all circuits $C$ with $C \subseteq F$. In other words,
$m_F$ is the least common multiple of those vertices of
$\,{\ini}_\prec(J_L) \,$ which lie on $F \cap H_\prec$.

The  labeled cell complex ${\ini}_\prec(\cHH_L)$ satisfies the
exactness criterion in [BS, Proposition 1.2]  because
the subcomplex of faces whose label
divides $ \xx^\aa \yy^\bb$ is contractible.
This follows from a result of Bj\"orner and Ziegler
[BLSWZ, Theorem 4.5.7]. Minimality is inherited from $\cHH_L$
since different faces in $\cHH_L$ have different labels.
\Box

\ref{initialRes} tells us that the minimal cellular resolution
$\,(\FF_{\cHH_L/L},\partial)\,$ of a unimodular Lawrence ideal $J_L$
is a {\it universal resolution\/} in the sense that it is stable with
respect to any term order $\prec$.  This generalizes the fact that the
minimal generators of $J_L$ are a {\it universal Gr\"obner basis\/}
(\ref{lawrencegens}).  In other words, the universal Gr\"obner basis
property extends from the generators of $J_L$ to all the higher 
syzygies.

\section{examples} Lawrence ideals arising from graphs

\noindent  Unimodular lattices arise naturally in the study of directed graphs
and matroids.  See Chapters 1 and 6 in [Whi1] for a first introduction.
Detailed information on unimodularity appears in
[Whi2, Chapter 3].  For instance, a famous theorem of Seymour 
[Whi2, Theorem 3.1.1(9)] states that every unimodular lattice can be
built up in a simple way (by duality, $1$-sums, $2$-sums and $3$-sums)
from graphic lattices and a certain five-dimensional lattice in
$\ZZ^{10}$.  In this section we discuss the minimal free resolutions
of cographic and graphic lattice ideals in combinatorial terms.

Let $G=(V,E)$ be a finite directed graph on $d$ vertices, which we
assume are labeled with numbers $1,2,\ldots, d$.  Suppose that $E$ has
$n$ edges.  The edge-node incidence matrix of $G$ is the $n \times
d$-matrix whose rows are $\, \ee_i-\ee_j\in \ZZ^d\,$ for $\,(i,j) \in
E$.  The image of this matrix in $\ZZ^n$ is denoted by $\, L_G$, and
is called the {\it graphic lattice\/} of $G$.  The orthogonal
complement of $L_G$ in $\ZZ^n$ is denoted by $\, L^*_G$, and is called
the {\it cographic lattice\/} of $G$.  In the language of matroid
theory, the graphic lattice $L_G$ is spanned by the {\it cocircuits}
of $G$, and the cographic lattice $L^*_G$ is spanned by the
{\it circuits\/} of $G$. The following classical result appears in
[Whi2, Theorem 1.5.3].

\proposition{graphsareuni}
The graphic and cographic lattices $L_G$ and $L^*_G$
are unimodular.

We write $J_G = J_{L_G}$ and $J^*_G = J_{L^*_G}$ for the two
unimodular Lawrence ideals associated with a directed graph $G$. We
call $J_G$ the {\it graphic ideal\/} and $J^*_G$ the {\it cographic
ideal\/}. In fact, these ideals depend only on the undirected graph
underlying $G$, and they can be described as follows.  Replace each
edge $(i,j)$ of the directed graph $G = (V,E)$ by two directed edges,
one from vertex $i$ to vertex $j$ and the other one from vertex $j$ to
vertex $i$. We associate the variables $x_{ij}$ and $x_{ji}$ with
these directed edges.  This gives a polynomial ring $S$ with $2n$
variables over $k$. We interpret the variable $x_{ji}$ as the
homogenizing variable for $x_{ij}$.  The pair of
variables $\{ x_{ij}, x_{ji} \}$ will play the same role as the pair
$\,\{ x_i,y_i \}\,$ in the previous sections.

We first discuss the graphic ideal of $G$. It can be computed
as an ideal quotient:
$$
J_G \quad = \quad
\biggl\langle
\prod_{j: (i,j) \in E} x_{ij} \,\,\, - \! \prod_{j: (i,j) \in E} x_{ji}
\,\,\, \, \big| \,\,\,\, i=1,2,\ldots,d  \,
 \biggr\rangle \,: \,
\langle \prod_{(r,s) \in E} x_{rs} x_{sr} \, \rangle^\infty.
$$
The $d$ binomials listed above correspond to a lattice basis of
$L_G$. Hence $J_G$ has codimension $d$.  The minimal generators of
$J_G$ are the {\it cocircuits} of the graph $G$.

\example{k5}
Consider the complete graph on five nodes, $\, G =  K_5$. 
The graphic ideal $J_{K_5}$ is generated by
five quartics such as
$\, x_{12}x_{13}x_{14}x_{15}-x_{21}x_{31}x_{41}x_{51}\,$
and ten sextics such as $\,
x_{13}x_{14}x_{15}x_{23}x_{24}x_{25}-x_{31}x_{32}x_{41}x_{42}x_{51}x_{52} $.
They correspond to the $15$ cocircuits of $K_5$.
The minimal free resolution of $J_{K_5}$
is given by a $4$-dimensional simplicial complex with $24$ facets,
$60$ tetrahedra, $50$ triangles and $15$ edges.

\vskip .2cm

We generalize this example by describing the minimal free
resolution of $J_{K_d}$, the graphic ideal of the complete
graph on $d$ nodes. Our construction is related to 
Lie algebra cohomology; a geometric version was found
independently by Bj\"orner and Wachs [BW].
Let $S$ be the polynomial ring over $k$ in the
$d(d-1)$ variables $x_{ij}$.  Let
$F_{r-1}$ denote the free $S$-module whose basis elements correspond
to ordered partitions $(A_1\!\mid\! A_2\!\mid\! \ldots\!\mid\! A_r)$ of the 
set $\set{1,\ldots, d}$, such that $1\in A_1$, for all $1\le r\le d-1$.

\theorem{barres}
The minimal resolution of the graphic ideal $J_{K_d}$ is the exact complex
$$ \FF_\bullet: \quad 0 \rightarrow F_{d-2}
\,\, {\buildrel \partial_{d-2} \over \longrightarrow}\,\, F_{d-3} \,\,
{\buildrel \partial_{d-3} \over \longrightarrow}\,\, F_{d-4} \,\,
\longrightarrow \ldots \longrightarrow \,\, F_2\,\,
 {\buildrel \partial_2 \over \longrightarrow}\,\, F_1 \,\,
 {\buildrel \partial_1 \over \longrightarrow}\,\, F_0 \rightarrow 0,
\eqno (5.1) $$
where the  differential $\partial_{r-1}$ acts on the
basis elements of $F_{r-1}$ by the ``cyclic rule''
$$\eqalign{
\partial_{r-1}&  (A_1\mid  A_2\mid \cdots\mid A_r)  \quad = \quad
 {(-1)}^{r+1}\prod_{i\in A_{r}}\prod_{j\in A_1}x_{ij} \cdot
(A_r\cup A_1\mid A_2 \mid \cdots\mid A_{r-1})\cr
 & +
\sum_{s=2}^r{(-1)}^s \! \prod_{i\in A_{s-1}}\prod_{j\in A_s}
x_{ij}\cdot (A_1\mid A_2\mid \cdots\mid A_{s-2}\mid A_{s-1}\cup A_s\mid
A_{s+1}\mid\cdots\mid A_r). \cr} $$

\proof
The central hyperplane arrangement $\cHH_{L_{K_d}}^0$
can be identified with the familiar {\it braid arrangement}
which consists of  the hyperplanes
$z_i = z_j$ in the $(d-1)$-dimensional vector space
$\setdef{(z_1,z_2\ldots,z_d) \in \RR^d}
{z_1 + \cdots + z_d = 0}$.
The $(r-2)$-faces of $\cHH_{L_{K_d}}^0$ are naturally labeled by the
ordered partition $(A_1\mid A_2\mid \ldots\mid A_r)$
of $\{1,2,\ldots,d\}$. More precisely, the
unique face having  a given point in its relative interior
is determined by  sorting the coordinates of that point.

We shall apply the initial ideal construction of \ref{initialRes}.
Let  $H_\prec$ denote the hyperplane $\{z_1 = 1 \}$ which represents
a lexicographic order with the variables $x_{12}, \ldots, x_{1d}$
being highest.  The faces of $\cHH_{L_{K_d}}^0$
which have bounded intersection with $H_\prec$ are
indexed  by those ordered partitions
 $(A_1\mid A_2\mid \ldots\mid A_r)$
which satisfy $ 1 \in A_1$.  The   simplicial complex
$\ini_\prec(\cHH_{L_{K_d}})$ is the
first barycentric subdivision of the $(d-2)$-simplex,
as shown for $d=4$ in \nextfig. The vertices of
$\ini_\prec(\cHH_{L_{K_d}})$  are indexed by the  monomials
$\, \prod_{i\in A_1}\prod_{j\in A_2}x_{ij} \,$
which represent partitions $(A_1,A_2)$ with $1 \in A_1$.
Higher dimensional faces are labeled by ordered partitions
with three or more parts.
The cellular resolution given  by the labeled
simplicial complex $\ini_\prec(\cHH_{L_{K_d}})$
has the format (5.1) with the differential mapping
$\, (A_1\mid  A_2\mid \cdots\mid A_r) \,$ to
$$ \sum_{s=2}^r{(-1)}^s \! \prod_{i\in A_{s-1}}\prod_{j\in A_s}
x_{ij}\cdot (A_1\mid A_2\mid \cdots\mid A_{s-2}\mid A_{s-1}\cup A_s\mid
A_{s+1}\mid\cdots\mid A_r).  $$
This complex is the minimal free resolution of the monomial ideal
$\ini_\prec(J_{K_d})$ as in \ref{initialRes}.
We lift this complex to the minimal free resolution
of $J_{K_d}$ by a construction as in  [PS, Theorem 5.4].
 Because $\cHH_{L_{K_d}}$ is simplicial, lifting amounts to
adding one more term to each differential, and that term is exactly the
remaining cyclic term, which is the first one listed in the statement of
\ref{barres}.
\Box

\capdraw{100}{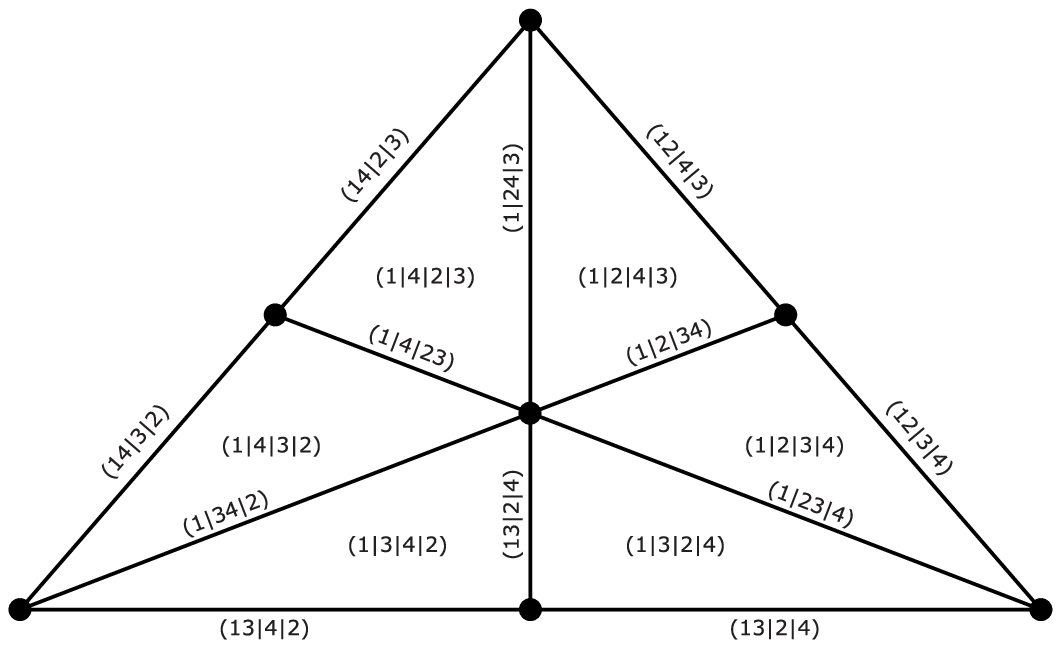}{
\stext{.5}{1}{\spad{3pt}{\small$x_{13}^{}x_{23}^{}x_{43}^{}$}}
\setext{.245}{.52}{\sepad{3pt}{3pt}{\small$x_{12}^{}x_{13}^{}x_{42}^{}x_{43}^{}$
}}
\swtext{.75}{.52}{\swpad{3pt}{3pt}{\small$x_{13}^{}x_{14}^{}x_{23}^{}x_{24}^{}$}
}
\wtext{.51}{.36}{\wpad{9pt}{\small$x_{12}^{}x_{13}^{}x_{14}^{}$}}
\ntext{0}{.025}{\npad{3pt}{\small$x_{12}^{}x_{32}^{}x_{42}^{}$}}
\ntext{.5}{.025}{\npad{3pt}{\small$x_{12}^{}x_{14}^{}x_{32}^{}x_{34}^{}$}}
\ntext{1}{.025}{\npad{3pt}{\small$x_{14}^{}x_{24}^{}x_{34}^{} $}}
}{An initial ideal of the complete graphic Lawrence ideal.}

We leave it to a future project to compute the
Betti numbers of the graphic ideal
$J_{G}$ for arbitrary directed graphs $G$. We expect
that there exist nice formulas in terms of the
characteristic polynomials of graphic arrangements
following [Za].

\vskip .1cm

Another open combinatorial problem is to find a 
formula for the Betti numbers of the {\it cographic ideals} $\, J^*_G $
which can be described as follows. Let $C = (C^+, C^-)$
be a signed circuit of the directed graph $G$;
see [BLSWZ, \S 1.1]. This means
that $C^+$ and $C^-$ are disjoint subsets of the set of edges
such that the edges in $C^+$ together with the reversals of the edges
in $C^-$ form a directed cycle which meets each vertex of $G$
at most once. Every signed circuit $C = (C^+,C^-)$ is coded into a binomial:
$$
\prod_{(i,j) \in C^+} \prod_{(k,l) \in C^-} x_{ij} x_{lk} \,\quad - \,\,
\prod_{(i,j) \in C^+} \prod_{(k,l) \in C^-} x_{ji} x_{kl} .$$
The cographic ideal $J^*_G$ is minimally generated by
these binomials where $C = (C^+,C^-)$ runs over all
signed circuits of the directed graph $G$.
For instance, if $G = K_4$ is the complete graph
on $\{1,2,3,4\}$ with edges $(i,j)$ for $i < j$, then
there are seven circuits  [BLSWZ, bottom of page 3] and
our  cographic ideal equals
$$
\eqalign{
J_{K_4}^*  \,\, =\,\,
\bigl\langle \, & \,
x_{12} x_{23} x_{31} - x_{21} x_{32} x_{13}  \,,\,\, 
x_{12} x_{24} x_{41} - x_{21} x_{42} x_{14} \,,\,\,
x_{13} x_{34} x_{41} - x_{31} x_{43} x_{14} ,\cr
&x_{23} x_{34} x_{42} - x_{32} x_{43} x_{24} \,,\,\,\,\,
 x_{12} x_{23} x_{34} x_{41} - x_{21} x_{32} x_{43} x_{14} \,,\,\cr
& x_{13} x_{32} x_{24} x_{41} - x_{31} x_{23} x_{42} x_{14} \,,\,\,\,\,
x_{12} x_{24} x_{43} x_{31} - x_{21} x_{42} x_{34} x_{13} \,\bigr\rangle \cr}
$$
Since the matroid of $K_4$ is self-dual, that is,
$K_4$ is self-dual under planar duality of graphs,
it follows that the cographic ideal $J_{K_4}^* $ is isomorphic
to the graphic ideal $J_{K_4} $ whose minimal resolution
is depicted in \ref{fig:initial}. For $d \geq 5$, however, the cographic ideal
$J_{K_d}^*$ is different from and
more complicated than  the graphic ideal $J_{K_d}$.

A special case of interest is the  complete bipartite
graph $K_{d,e}$. Its cographic  ideal $J^*_{K_{d,e}}$
is the Lawrence liftings of
the ideal of $2 \times 2$-minors of a generic matrix.  While the
minimal free resolution of the ideals of $2 \times 2$-minors 
depends on the characteristic of the
base field, a result of Hashimoto [RW], their Lawrence liftings
have minimal free resolutions  which are cellular and
characteristic free, by \ref{ZeroOneagain}.

\section{diagonal} The diagonal embedding of a unimodular toric variety

\noindent 
We  next present a different geometric interpretation
of the Lawrence ideal $J_L$. It concerns the diagonal embeddings of a
toric variety, written in the homogeneous coordinates of Audin-Cox
[Au],[Cox]. This connection first appeared in [Oda].
As an application,  Beilinson's spectral sequence is presented for
unimodular toric varieties.

Let $\Sigma$ be a complete fan in  the lattice $\ZZ^m$
and $X$ the associated complete normal toric variety; see e.g.~[Fu].
Let $\bb_1,\dots,\bb_n \in \ZZ^m$ be the primitive generators
of the one-dimensional cones of $\Sigma$,
and let $B$ be the $n \times m$-matrix with row vectors $\bb_i$.
Each $\bb_i$ determines a torus-invariant Weil divisor $D_i$ on $X$,
and the group $\Cl(X)$ of torus-invariant Weil divisors modulo
linear equivalence has the presentation
$$
0 \,\, \longrightarrow \,\,\ZZ^m
\,\,\, {\buildrel B \over \longrightarrow }  \,\,\, \ZZ^n
\,\,\, {\buildrel \pi \over \longrightarrow }  \,\,\, \Cl(X)
\,\,\, \longrightarrow 0,
$$
%
where $\pi$ takes the $i$-th standard basis vector
of $\ZZ^n$ to  the linear equivalence 
class $[D_i]$ of the corresponding divisor $D_i$.
If the divisor class group is torsion free, so that
$\Cl(X) = \ZZ^{n-m}$, then we express $\pi$ by
a matrix $A$, and we are back in (2.1).

Audin [Au] and Cox [Cox] define the {\it homogeneous coordinate ring\/}
of $X$ to be the polynomial ring $R=k[x_1,\ldots,x_n]$, equipped
with a grading by the abelian group $\Cl(X)$ via the morphism $\pi$
above. Hence a monomial $\xx^\aa\in R$ has degree
$[\sum_{i=1}^n a_i D_i]\in\Cl(X)$. The fan $\Sigma$ is encoded
in the {\it irrelevant ideal} $B_\Sigma \subset R$ 
whose monomial generators correspond to complements of
facets of $\Sigma$. Coherent sheaves on the toric variety $X$ are 
represented by  $B_\Sigma$-torsion-free
$\Cl(X)$-graded modules over $R$; see [Cox] for the
case where $\Sigma$ is simplicial, and [Mu] for the general case.
Closed subschemes of $X$ are defined by $B_\Sigma$-saturated
$\Cl(X)$-graded ideals of $R$.  Furthermore, for any
$\Cl(X)$-graded $R$-module $T$, and any $\aa\in \Cl(X)$, we define its
{\it twist} $T(\aa)$ to be the graded $R$-module with components
$T(\aa)_\bb:=T_{\aa+\bb}$. Let  $\cOO_X(\aa)$ denote the coherent
sheaf on $X$ corresponding to the twisted $R$-module $R(\aa)$.

The toric variety $X \times X$ has the 
homogeneous coordinate ring
$$  S \quad = \quad R \, \otimes_k \, R 
\quad = \quad k[x_1,\ldots,x_n, \, y_1,\ldots,y_n]  $$
The diagonal embedding $X\subset X\times X$ defines
a closed subscheme, and is  represented by
a $\Cl(X)\times\Cl(X)$-graded ideal $I_X$ in $S$.
That ideal  $I_X $ is the kernel of  the
 natural map $\,\, S =  R \otimes R \, \rightarrow \,
k[\Cl(X)]\otimes R \, ,\,\,
\xx^\uu \yy^\vv = 
\xx^\uu \otimes \xx^\vv \, \mapsto \,
[\uu]\otimes\xx^{\uu+\vv}$. Explicitly,
$$I_X \quad =  \quad \idealdef{\xx^\uu \yy^\vv - \xx^\vv \yy^\uu}
{\pi(\uu)=\pi(\vv)\, \, \hbox{in} \,\, \Cl(X)} \quad \subset \quad S .$$
This formula proves the following observation.

\proposition{LawrenceGivesDiagonal}
The ideal $I_X \subset  S $ defining the diagonal embedding
$X \subset X \times X$ equals the Lawrence ideal
$J_L$ for the lattice $L=\im(B)=\ker(\pi)$
of principal divisors.

The basic example is projective space $X = \PP^m$.
Here $n= m+1$, $\, L = \ker (\, 1 \, \,1 \, \,1 \, \, \cdots \, 1 \,) $,
$\Cl(X) = \ZZ^{1}$, and $I_X = J_L$ is the
ideal (1.1) of $2 \times 2$-minors of a $\, 2 \times n$-matrix
of indeterminates. The minimal free resolution of  $I_X = J_L$ is 
an Eagon-Northcott complex.
Example 10 in [Oda] demonstrates how
Beilinson's spectral sequence for $\PP^m$
can be derived from this observation.

The results in this paper provide a cellular model for this
Eagon-Northcott complex as a hypersimplicial complex (Figure 1).
We will show that this derivation of Beilinson's spectral sequence extends 
to the general setting of \ref{unimod-toric-syz}. Thus our cellular
resolutions provide a partial solution to [Oda, Problem 6].

We say that a toric variety $X$ is {\it unimodular}
if its lattice $L$ of principal divisors is a unimodular
sublattice of $\ZZ^n$. This condition is equivalent,
by [Stu, Remark 8.10], to  saying that the toric variety 
$X$ is smooth and any other variety obtained by toric flips and 
flops is smooth as well. Toric varieties with this property
appear frequently in representation theory and in integer programming.
For example, the toric varieties associated with
{\it transportation polytopes} and {\it products of minors}
are unimodular. These toric varieties were studied recently 
by Babson and Billera [BB].

Suppose now that $X$ is a unimodular toric variety.
Then \ref{unimod-toric-syz} constructs an $m$-dimensional cellular model
for the $\Cl(X)\times\Cl(X)$-graded minimal 
resolution for the ideal 
$\, I_X \subset S \,$ of the diagonal embedding
of $X$ in $X\times X$:
$$L_\bullet: \quad 0\longrightarrow L_m \longrightarrow \ldots \longrightarrow
L_2\longrightarrow L_1 \longrightarrow L_0=\cOO_{X\times X}
\longrightarrow \cOO_{X\times X}/\cII_X \longrightarrow 0. $$
Here we identify $I_X$ with the associated coherent sheaf
$\cII_X \subset\cOO_{X\times X}$, which is the ideal sheaf of
the diagonal embedding of $X$. Note that $S/I_X$ is Cohen-Macaulay.

Since $\cII_X$ is a $\Cl(X)\times\Cl(X)$-homogeneous ideal, we can write
$\, L_i \, = \,\oplus_{j=1}^{m_i}\cOO(-\aa_{ij},-\bb_{ij})\,$
with $\aa_{ij}, \bb_{ij}\in\Cl(X)$. Here
$\cOO(-\aa_{ij},-\bb_{ij})=\cOO(-\aa_{ij})\boxtimes\cOO(-\bb_{ij})=
p_1^*(\cOO(-\aa_{ij}))\otimes p_2^*(\cOO(-\bb_{ij}))$
denotes the exterior tensor product on $X\times X$ of 
$\cOO(-\aa_{ij})$ and $\cOO(-\bb_{ij})$, and 
$p_i$ are the projections on the factors of $X\times X$.

Let $\cFF$ be any coherent sheaf on $X$. By tensoring the above free
resolution $L_\bullet$ with $p_1^*(\cFF)$ we obtain a resolution
$C_\bullet$ for the restriction of $p_1^*(\cFF)$ to the diagonal
$X \subset X\times X$. As in the case of projective space [Bei], pushing
via ${p_2}_*$ a Cartan-Eilenberg injective resolution of $C_\bullet$
yields a double complex of $\cOO_X$-modules, whose total chain complex
has as cohomology the hyperdirect image ${\mathbf
R}{p_2}_*(C_\bullet)$.  One of the spectral sequences belonging to the
double complex ${\mathbf R}{p_2}_*(C_\bullet)$ yields
the following version of the
{\it Beilinson Spectral Sequence}. It
says that any coherent sheaf $\cFF$ can be reconstructed from the
knowledge of the cohomology of a few of its $\Cl(X)$-twists.

\theorem{beilinson.ss} Let $\cFF$ be a coherent sheaf on a complete 
unimodular toric variety $X$. Then there is a 
(third quadrant) spectral sequence with $E_1$ term
$$E_1^{pq} \quad = \quad \bigoplus_j
H^q(\cFF\otimes\cOO_X(-\aa_{pj}))\otimes\cOO_X(-\bb_{pj})$$
which converges to the associated graded sheaf of 
a filtration of $\cFF$.
\Box

In the case of $\PP^m$, Beilinson's spectral sequence induces an
equivalence between the derived category $\DD^b(\Coh(\PP^m))$ of
bounded complexes of sheaves with coherent cohomology, and the derived
category of modules over the exterior algebra in $m+1$ generators.
An analogous simple description of $\DD^b(\Coh(X))$ does not hold over
an arbitrary unimodular toric variety $X$, since the bundle
$\cEE=\oplus_{i,j}\cOO(-\aa_{ij})$ defined as the direct sum (without
repetitions) of all ``left'' summands in the above resolution of the
ideal sheaf of the diagonal is in general not exceptional, i.e.
$\Ext^i_{\cOO_X}(\cEE,\cEE)\ne 0$ for some $i>0$. 
See [AH] for a related approach to  describing $\DD^b(\Coh(X))$.

\references

\itemitem{[AH]} K.~Altmann, L.~Hille: 
Strong exceptional sequences provided by quivers.
{\sl Algebr. Represent. Theory} {\bf 2} (1999) 1-17.

\itemitem{[Au]} M.~Audin: {\sl The Topology of Torus Actions on
Symplectic Manifolds}, Progress in Math. {\bf 93},
Birkh\"auser, 1991.

\itemitem{[BB]} E.~Babson, L.J.~Billera:
The geometry of products of minors,
{\sl Discrete Comput. Geometry} {\bf 20} (1998) 231--249

\itemitem{[BPS]} D.~Bayer, I.~Peeva, B.~Sturmfels:
Monomial resolutions, {\sl Math. Res. Lett.} {\bf 5} (1998), no. 1-2,
31--46.

\itemitem{[BS]} D.~Bayer, B.~Sturmfels:
Cellular resolutions of monomial modules,
{\sl J. Reine Angew. Math.} {\bf 502} (1998),
123--140.

\itemitem{[Bei]} A.~Beilinson:
Coherent sheaves on $\PP^N$ and problems of linear algebra,
{\sl Funct. Anal. Appl.}, {\bf 12}, (1978), 214--216.

\itemitem{[BLSWZ]} A.~Bj\"orner, M.~Las~Vergnas,
B.~Sturmfels, N.~White and G.~Ziegler:
{\sl Oriented Matroids}, Cambridge University Press, 1993.

\itemitem{[BW]} A.~Bj\"orner, M.~Wachs:
Geometrically constructed bases
for intersection lattices of real arrangements, preprint 1999.

\itemitem{[Cox]} D.A.~Cox: The homogeneous coordinate ring of a
toric variety,
{\sl J. Algebraic Geom.} {\bf 4} (1995), no. 1,
17--50.

\itemitem{[Fu]} W.~Fulton:
{\sl Introduction to Toric Varieties}, Princeton
Univ.~Press, 1993.


\itemitem{[Mu]} M.~Musta\c t\u a:
Vanishing theorems on toric varieties, preprint 1999.

\itemitem{[Oda]} T.~Oda: Problems on Minkowski sums of
convex lattice polytopes, Manuscript, 1997.

\itemitem{[PS]} I.~Peeva, B.~Sturmfels: 
Generic lattice ideals,
{\sl J. Amer. Math. Soc.}, {\bf 11}, (1998), no. 2,
363--373.

\itemitem{[RW]} J.~Roberts, J.~Weyman: 
A short proof of a theorem of M. Hashimoto. 
{\sl J. Algebra} {\bf 134} (1990), no.~1, 144--156.  

\itemitem{[Stu]} B.~Sturmfels:
{\sl Gr\"obner Bases and Convex Polytopes},
AMS University Lecture Series,  Vol. {\bf 8}, Providence RI,
1995.

\itemitem{[Whi1]} N.~White: {\sl Theory of Matroids},
Encyclopedia of Mathematics and its Applications, {\bf 26},
Cambridge University Press, 1986.

\itemitem{[Whi2]}  N.~White: {\sl Combinatorial Geometries},
Encyclopedia of Mathematics and its Applications, {\bf 29},
Cambridge University Press, 1987.

\itemitem{[Za]} T.~Zaslavsky: Facing up to arrangements: face-count formulas
for partitions of space by hyperplanes, {\sl Mem. Amer. Math. Soc.} {\bf 1}
(1975), no.~154.

\itemitem{[Zi]} G.~Ziegler: {\sl Lectures on Polytopes},
Graduate Text in Mathematics {\bf 152}, Springer Verlag, 1997.

\vskip .6cm
\widow{.1}

\noindent
Dave Bayer, Department of Mathematics,
Barnard College,
Columbia University,
New York, NY 10027, USA,
{\tt
bayer@math.columbia.edu}.

\vskip .2cm

\noindent
Sorin Popescu, Department
of Mathematics,
Columbia University, New York, NY 10027, USA, and
Department of Mathematics, SUNY at Stony Brook, Stony Brook, NY 11794,
{\tt
psorin@math.columbia.edu}.

\vskip .2cm

\noindent
Bernd
Sturmfels,
Department of Mathematics,
University of California,
Berkeley,
CA 94720, USA,
{\tt bernd@math.berkeley.edu}.
\bye